# MODELS AND METHODS FOR THREE EXTERNAL BALLISTICS INVERSE PROBLEMS

Natalya Arutyunova*, Aidar Dulliev*, Vladislav Zabotin*

October 5, 2016

**Abstract**

We consider three problems of selecting optimal gun barrel direction (or those of selecting optimal semi-axis position) when firing an unguided artillery projectile on the assumption that the gun barrel semi-axis can move in a connected nonconvex cone having a non-smooth lateral surface and modelling visibility zone restrictions. In the first problem, the target is in the true horizon plane of the gun, the second and the third problems deal with some region of 3D space. A distinctive feature of the models is that the objective functions are ε-Lipschitz ones. We have constructed a unified numerical method to solve these problems based on an algorithm of projecting a point onto ε-Lipschitz level function set. A computer program has been based on it. A series of numerical experiments on each problem has been carried out. App. 1.

Russian text – pp. 1–26, English text – pp. 27–37.

* Department of Applied Mathematics and Informatics, Kazan National Research Technical University named after A.N. Tupolev-KAI, Russian Federation. E-mail: {NKArutyunova@kai.ru, AMDulliev@kai.ru, VIZabotin@kai.ru}.

# МОДЕЛИ И МЕТОДЫ ДЛЯ ТРЕХ ОБРАТНЫХ ЗАДАЧ ВНЕШНЕЙ БАЛЛИСТИКИ

Н.К. Арутюнова*, А.М. Дуллиев*, В.И. Заботин*

**Аннотация**

Рассматриваются три математические модели задачи выбора оптимального направления ствола орудия при стрельбе неуправляемым снарядом в предположении, что полуось ствола может перемещаться в связном невыпуклом конусе, имеющем негладкую боковую поверхность и моделирующем ограничения на зону видимости. В первой задаче цель расположена в плоскости истинного горизонта орудия, во второй и третьей – в некоторой области пространства. Отличительной особенностью моделей является ε-липшицевость целевых функций. Построен единый численный метод решения поставленных задач, базирующийся на одном алгоритме проектирования точки на множество уровня ε-липшицевой функции. На его основе составлена программа для ЭВМ. По каждой из задач проведена серия вычислительных экспериментов. Библ. 6.

## 1. ВВЕДЕНИЕ

Задача стрельбы неуправляемым снарядом, в которой требуется минимизировать расстояние от точки падения снаряда до цели, в настоящее время изучена в теории внешней баллистики достаточно хорошо при условии пренебрежения атмосферой и при условии того, что перемещение ствола орудия ограничено лишь плоскостью местного горизонта, а поверхность Земли считается сферической [1].

Однако в реальной ситуации положение полуоси ствола орудия можно менять произвольно, как правило, лишь в пределах некоторого связного невыпуклого конуса, имеющего негладкую боковую поверхность и возникающего в условиях, сужающих выбор направления ствола из-за наличия каких-либо препятствий.

* Кафедра прикладной математики и информатики, Казанский национально-исследовательский технический университет им. А.Н. Туполева, Российская Федерация. E-mail: {NKArutyunova@kai.ru, AMDulliev@kai.ru, VIZabotin@kai.ru}.

Кроме ограничений на перемещение ствола орудия, нередко приходится учитывать тот факт, что цель может находиться вне плоскости истинного горизонта орудия и, возможно, на некоторой поверхности, определяющей рельеф местности на заданной территории. В последнем случае проблема выбора положения полуоси ствола становится гораздо сложнее. Связано это главным образом с тем, что минимальное в евклидовой метрике расстояние от цели до точки падения снаряда не всегда соответствует оптимальному или даже близкому к оптимальному выстрелу. Например, в тех случаях, когда цель и точка падения снаряда отделены друг от друга какой-либо существенной преградой.

Математически модель может быть описана следующим образом. Обозначим: dist – евклидово расстояние, $l$ – луч с началом в точке $O$ нахождения снаряда, $N$ – точка падения снаряда, $N = N(l)$, $M$ – цель, $K$ – конус с вершиной в точке $O$, $D'$ – множество трехмерного евклидова пространства, $M \in D'$. Требуется

$$\min_{l} \operatorname{dist}(N(l), M),$$

$$\text{subject to } l \in K,$$

$$N(l) \in D'.$$

Возможны постановки задачи и с другими целевыми функциями. Так, ниже рассмотрена задача, в которой в качестве целевой функции фигурирует расстояние между траекторией полёта снаряда и точкой $M$.

Сразу оговоримся, что в настоящей статье, мы не учитываем сопротивление атмосферы. Однако, предлагаемые и анализируемые нами модели, во-первых, могут использоваться для стрельб тяжелыми неуправляемыми снарядами. Во-вторых, при стрельбе лёгкими снарядами полученные решения могут считаться опорными и уточняться с помощью соответствующего аппарата.

Ниже ставятся и исследуются три задачи, для них строится единый алгоритм решения и на его основе составлены программы для ЭВМ. По каждой из задач проведена серия вычислительных экспериментов.

## 2. ПОСТАНОВКА ЗАДАЧ

Всюду ниже орудие будем моделировать материальной точкой $T$, а направление ствола орудия – лучом с началом в этой точке; несферичностью Земли будем пренебрегать.

В качестве системы координат возьмём декартову систему координат $Oxyz$, начало $O$ которой совпадает с точкой $T$, плоскость $Oxy$ – плоскость истинного горизонта в точке $O$, направления осей $Ox$ и $Oy$ будет уточнено ниже, ось $Oz$ направлена перпендикулярно $Oxy$ и противоположно вектору ускорения силы тяжести $\boldsymbol{g}$. Всюду в данной статье будем использовать следующие обозначения: $\|\cdot\|_p$ – $p$-норма в $\mathbb{R}^n$, $n \geq 1$, при $p = 2$ нижний индекс $p$ будем опускать; int, fr – операторы внутренности и границы в $\mathbb{R}^n$; $\mathrm{Pr}_i\, A$ – проекция множества $A$ на $i$-ю координатную ось (здесь $i$ равно 1, или 2, или 3 для осей $x$, $y$, $z$ соответственно); $\mathrm{Lip}(X; Y)$ – множество липшицевых на $X$ функций $f$: $X \to Y$ с константой Липшица $\mathrm{lip}(f)$.

Мы будем считать, что направление (положение полуоси) ствола орудия может свободно выбираться в пределах замкнутого конуса $K$, с вершиной в точке $O$ и не содержащего вертикальной образующей. Пренебрегая длиной ствола орудия, начальную скорость $\boldsymbol{v}$ снаряда будем считать приложенной в точке $O$ и постоянной по модулю. Другими словами, $\boldsymbol{v} \in S \cap K$, где $S = \{\boldsymbol{v} : \|\boldsymbol{v}\| = v_0\} \subset \mathbb{R}^3$ – сфера радиуса $v_0$ с центром в точке $O$. Через $\nu$ обозначим дробь $v_0 / \sqrt{\|\boldsymbol{g}\|}$.

Всюду далее выбор направления ствола орудия считаем эквивалентным выбору вектора $\boldsymbol{v} = (v_x, v_y, v_z)$ с заданным фиксированным значением $v_0$, который будем задавать посредством сферической системы координат $(v_0, \varphi, \psi)$ в виде

$$\begin{cases} v_x = v_0 \cos\psi \cos\varphi, \\ v_y = v_0 \cos\psi \sin\varphi, \\ v_z = v_0 \sin\psi; \end{cases} \tag{2.1}$$

$$(\varphi; \psi) \in E \subset \Theta := \{(\varphi; \psi) : \varphi \in (-\pi; \pi],\ \psi \in (-\pi/2; \pi/2)\},$$

где множество $E$ – замкнутое и соответствует множеству $S \bigcap K$ в координатах $(\varphi, \psi)$.

Влияние сопротивления воздуха на траекторию движения снаряда в данной статье не учитывается. Отсюда следует, что снаряд движется по параболической траектории в плоскости стрельбы.

Рассмотрим три задачи выбора оптимального направления ствола орудия при стрельбе неуправляемым снарядом по заданной цели. При этом будем считать, что цель $M$ задается точкой, координата $y$ которой равна нулю и которая принадлежит множеству, ограниченному параболоидом поражения [1], получающимся при отсутствии ограничений на перемещение ствола орудия. (Случай, когда цель может находиться вне пределов параболоида поражения, а также задача оптимального совместного выбора направления ствола орудия и скорости снаряда, в настоящей работе авторами не рассматриваются.)

Среди всех точек на траекториях снаряда будем рассматривать те, координата $x$ которых удовлетворяет неравенству $x \geq \kappa$, $\kappa = \mathrm{const}$, $\kappa > 0$. Практическое значение величины κ может быть различным, например, определяться радиусом поражения снаряда или же безопасным расстоянием от позиции стрельбы до поражаемого объекта (цели).

В зависимости от конкретной задачи другие ограничения на значения $x$ и $y$, а также ограничения на $z$ будут оказываться различными.

**Задача I.** Пусть цель, по которой ведётся стрельба, задаётся точкой $M$ в плоскости $xOy$ с координатами $(a; 0)$, $a > \kappa$. В предположении, что точка $N$ падения снаряда также расположена в плоскости $xOy$ и имеет в ней координаты $(x; y)$, а также отсутствуют препятствия для всех траекторий, определяемых вектором начальной скорости $\boldsymbol{v} \in S \bigcap K$, требуется среди этих траекторий выбрать такую, для которой евклидово расстояние между точками $M$ и $N$ минимально.

Расстояние $r$ от точки $O$ до точки $N$ $(x; y)$ в соответствии с [1] вычисляется по формуле:

$$r = \sqrt{x^2 + y^2} = v^2 \sin 2\psi. \tag{2.2}$$

Учитывая (2.1)–(2.2) получаем, что координаты $x$ и $y$ находятся из соотношений:

$$x = r\cos\varphi = v^2 \sin 2\psi \cos\varphi; \\ y = r\sin\varphi = v^2 \sin 2\psi \sin\varphi. \qquad (2.3)$$

Определяемое ими отображение $(\varphi, \psi) \mapsto (x, y)$ из $\Theta$ в $\mathbb{R}^2$ обозначим через $h_{\mathrm{I}}$.

Учитывая фиксированность модуля начальной скорости движения снаряда, из формулы (2) можно установить верхнее ограничение на дальность полёта $\rho_{\mathrm{I}} = v^2$. Введём обозначение

$$W_{\mathrm{I}} = \left\{ (x; y) : \sqrt{x^2 + y^2} \le \rho_{\mathrm{I}}; x \ge \kappa \right\}. \qquad (2.4)$$

Очевидно, что допустимое множество точек падения снаряда задаётся множеством $h_{\mathrm{I}}(E)$.

Таким образом, направление ствола орудия, задаваемое переменными $(\varphi; \psi)$ и соответствующее оптимальной траектории снаряда, находится из решения задачи:

$$\| M - N \|^2 \to \min, \\ \text{s.t.} \quad N \in h_{\mathrm{I}}(E) \cap W_{\mathrm{I}}. \qquad (2.5)$$

В настоящей статье множество $E$ будет описываться с помощью функциональных неравенств.

Пусть $[\theta_1; \theta_2] \subset (-\pi/2; \pi/2)$, а функции $g_1(\varphi)$ и $g_2(\varphi)$ таковы, что $-\pi/2 < g_1(\varphi) \le g_2(\varphi) < \pi/2$, $\varphi \in [\theta_1; \theta_2]$, положим

$$E = \left\{ (\varphi; \psi) \in \Theta : \varphi \in [\theta_1; \theta_2], g_1(\varphi) \le \psi \le g_2(\varphi) \right\},$$

или то же самое:

$$E = \left\{ (\varphi; \psi) \in \Theta : g(\varphi; \psi) \le 0 \right\},$$

где

$$g(\varphi; \psi) := \max\left\{ \theta_1 - \varphi, \varphi - \theta_2, g_1(\varphi) - \psi, \psi - g_2(\varphi) \right\}. \qquad (2.6)$$

Заметим, что задача I в определенном смысле допускает обобщение на случай, когда цель, по которой ведётся стрельба, состоит из множества

$\mathcal{M} = \{ M_1, M_2, \ldots, M_n \}$ точек $M_i = (a_i; b_i)$, $i = 1, \ldots, n$, $n > 1$, плоскости $xOy$. Действительно, при однократном выстреле в качестве цели $M$ можно принять, например, чебышёвский центр множества $\mathcal{M}$.

**Задача II.** Пусть цель, по которой ведётся стрельба, задаётся одной точкой $M$ $(a; 0; c) \in \mathbb{R}^3$. Предположим, что цель и все возможные точки падения снаряда расположены не ниже плоскости $z = z_{\min}$ ($z_{\min} < 0$). Как и в задаче I, будем считать, что $a > \kappa$. В зависимости от ограничений, накладываемых на точку $M$, выделим две подзадачи, называемые далее задача II.а и задача II.б.

**II.а.** Точка $M$ – произвольная точка (т.е. может быть расположена в воздухе), а выбор траектории полёта снаряда ограничивается лишь множеством $E$. Требуется из всех допустимых траекторий выбрать ту, для которой расстояние от неё до цели $M$ минимально по сравнению с расстояниями от остальных траекторий.

**II.б.** Точка $M$ лежит на границе fr $D$ некоторого множества $D$, задаваемого неравенством $H$ $(x; y; z) \le 0$, причём $H(0; 0; 0) \ge 0$ и функция $H$ непрерывна в $\mathbb{R}^3$. Эта функция определяет рельеф местности (т.е. цель может располагаться на поверхности земли/воды или некоторого неподвижного наземного объекта, а орудие находится на поверхности земли / воды либо в воздухе). Предполагается, что точка $N$ $(x; y; z)$ падения снаряда находится на fr $D$ и участок траектории снаряда от точки выстрела $O$ до точки $N$ не имеет общих точек с int $D$. Кроме того, будем считать, что из $z \le z_{\min}$ следует $H$ $(x; y; z) \le 0$. Требуется среди всех таких допустимых траекторий выбрать ту, для которой отрезок $MN$ не пересекается с int $D$ и имеет минимальную длину, или убедиться в отсутствии таковой. Содержательно это означает, что точка падения снаряда должна быть в прямой видимости от цели и располагаться при этом максимально близко к ней.

Цель $M$ может располагаться выше ($c > 0$) либо ниже (при $c < 0$) плоскости $xOy$, причем в последнем случае оптимальное значение для угла $\psi$ может оказаться меньше нуля.

Прежде, чем привести аналитические формулировки задач II.а и II.б, сделаем несколько важных замечаний.

В отличие от задачи I определение однозначного отображения из $\Theta$ в $\mathbb{R}^3$, аналогичного $h_{\text{I}}$, в задачах II.а и II.б либо вообще невозможно, либо весьма проблематично. Связано это с тем, что цель и область поражения в задачах II.а,б расположены в $\mathbb{R}^3$. А именно, в задаче II.а каждой точке $(\varphi;\ \psi) \in E$ соответствует не одна, а целое множество точек траектории полёта снаряда; что касается задачи II.б, аналитическое определение в ней координат точки $N\,(x;\ y;\ z)$ для заданных $(\varphi;\ \psi) \in E$, осуществимо, по-видимому, лишь в некоторых специальных случаях, определяемых видом множества $D$.

Оказывается, вместо многозначного отображения из $\Theta$ в $\mathbb{R}^3$ удобно рассматривать однозначное отображение $h_{\text{II}}\colon \Theta\times\mathbb{R}^{+}\to\mathbb{R}^{3}$, которое каждой точке $(\varphi;\psi;r)\in\Theta\times\mathbb{R}^{+}$ ставит в соответствие точку $(x;\ y;\ z) \in \mathbb{R}^3$, расположенную на отвечающей параметрам $(\varphi;\ \psi)$ траектории снаряда, при этом таким образом, что расстояние от точки $O$ до проекции точки $(x;\ y;\ z) \in \mathbb{R}^3$ на плоскость $xOy$ равно $r$.

Согласно теории внешней баллистики, такое отображение $h_{\text{II}}$ определяется равенством

$$h_{\text{II}}(\varphi,\psi,r)=\left\{(x,y,z): x=r\cos\varphi;\, y=r\sin\varphi;\, z=r\,\text{tg}\,\psi-\left(1+\text{tg}^2\psi\right)(2v^2)^{-1}r^2\right\}. \quad (2.7)$$

Несмотря на однозначность отображений $h_{\text{I}}$ и $h_{\text{II}}$, обратные к ним отображения $h_{\text{I}}^{-1}$ и $h_{\text{II}}^{-1}$ являются многозначными: заданной точки падения снаряда можно достичь путём выбора одного из двух подходящих значений угла наклона ствола орудия $\psi\in(-\pi/2;\pi/2)$. Впрочем, величина угла $\varphi$ при этом определяется однозначно. Отображения $h_{\text{I}}^{-1}$ и $h_{\text{II}}^{-1}$ в дальнейшем играют важную роль в построении функциональных ограничений для задач I и II.а,б. Причина, по которой данное построение оказывается возможным, заключается в том, что $h_{\text{I}}^{-1}$ и $h_{\text{II}}^{-1}$ допускают весьма простую параметризацию однозначными отображениями. А именно, из соотношений (2.3), (2.7) вытекает, что существуют однозначные отображения $h_{\mathrm{I}\,j}^{-1}$, $h_{\mathrm{II}\,j}^{-1}$, $j \in \{1;\, 2\}$ такие, что:

$$h_k^{-1} = h_{k1}^{-1} \bigcup h_{k2}^{-1}, \quad \mathrm{dom}\left(h_{kj}^{-1}\right) = \mathrm{dom}\left(h_k^{-1}\right), \quad k \in \{\mathrm{I};\mathrm{II}\}, \quad j \in \{1;2\},$$
$$\mathrm{Pr}_2\,\mathrm{Im}\left(h_{\mathrm{I}1}^{-1}\right) = \left(0;\pi/4\right], \qquad \mathrm{Pr}_2\,\mathrm{Im}\left(h_{\mathrm{I}2}^{-1}\right) = \left[\pi/4;\pi/2\right), \tag{2.8}$$
$$\mathrm{Pr}_2\,\mathrm{Im}\left(h_{\mathrm{II}1}^{-1}\right) = \left(-\pi/2;\pi/2\right), \qquad \mathrm{Pr}_2\,\mathrm{Im}\left(h_{\mathrm{II}2}^{-1}\right) = \left(0;\pi/2\right).$$

Явные формулы для вычисления $h_{kj}^{-1}$ будут приведены позже.

Поскольку все точки падения снаряда расположены выше плоскости $z = z_{\min}$, то в соответствии с (2.7) значение $r$ может быть ограничено сверху величиной

$$\rho_{\mathrm{II}} = v\sqrt{v^2 - 2z_{\min}}.$$

С учётом этого, а также огибающей семейства возможных траекторий введём множество:

$$W_{\mathrm{II}} = \left\{ (x;y;z) : \sqrt{x^2 + y^2} \le \rho_{\mathrm{II}};\, x \ge \kappa,\, z_{\min} \le z \le \frac{1}{2}\left( v^2 - \frac{x^2 + y^2}{v^2} \right) \right\}, \tag{2.9}$$

которое представляет собой множество достижимых снарядом точек.

Сделанные замечания позволяют представить аналитическую формулировку задачи II.а в виде

$$\begin{gathered} \left\| M - N \right\|^2 \to \min, \\ \text{s.t. } N \in h_{\mathrm{II}}\left(E\right) \bigcap W_{\mathrm{II}}. \end{gathered} \tag{2.10}$$

Для задачи II.б имеем следующую формулировку

$$\begin{aligned} & \left\| M - N \right\|^2 \to \min, \\ \text{s.t.}\quad & N \in h_{\mathrm{II}}\left(E\right) \bigcap W_{\mathrm{II}} \bigcap \mathrm{fr}\, D, \\ & \exists\, j \in \{1;2\} : h_{\mathrm{II}}\left(\mathrm{Pr}_1 h_{\mathrm{II}\,j}^{-1}(N); \mathrm{Pr}_2 h_{\mathrm{II}\,j}^{-1}(N); [0; \mathrm{Pr}_3 h_{\mathrm{II}\,j}^{-1}(N)]\right) \bigcap \mathrm{int}\, D = \varnothing, \\ & \left\{ \lambda N + (\lambda - 1) M : \lambda \in [0;1] \right\} \bigcap \mathrm{int}\, D = \varnothing. \end{aligned} \tag{2.11}$$

Согласно постановке задачи II.б последние два равенства в (2.11) означают, соответственно, что участок траектории снаряда от точки выстрела $O$ до точки $N$ и отрезок $MN$ не пересекаются с $\mathrm{int}\, D$.

Обозначим через $N_*$ решение задачи (2.11).

К сожалению, постановка задачи (2.11) допускает существование общих, помимо $N_*$, точек у участка траектории снаряда от точки $O$ до точки $N_*$ или отрезка $MN_*$ с множеством fr $D$. Конечно, к ограничениям задачи (2.11) можно добавить соответствующие условия, исключающие подобную ситуацию. Однако, как несложно показать, их формализация не позволяет в дальнейшем корректно пользоваться каким-либо методом минимизации, поскольку нарушается замкнутость множества ограничений, а попытка построить его замыкание вновь приводит к задаче (2.11). В качестве рецепта решения данной проблемы можно предложить следующую процедуру: решить задачу (2.11) и если вдруг окажется, что для ее решения $N_*$ указанное в начале абзаца условие выполнено, то провести поиск другого решения задачи (2.11) в объединении достаточно малых замкнутых «кольцевых» окрестностях точек $\left(\mathrm{Pr}_1 h^{-1}_{\mathrm{II},1}(N_*); \mathrm{Pr}_2 h^{-1}_{\mathrm{II},1}(N_*)\right)$ и $\left(\mathrm{Pr}_1 h^{-1}_{\mathrm{II},2}(N_*); \mathrm{Pr}_2 h^{-1}_{\mathrm{II},2}(N_*)\right)$, соответственно.

## 3. МИНИМИЗАЦИЯ ЦЕЛЕВЫХ ФУНКЦИЙ

Поставленные задачи (2.5), (2.10), (2.11) в общем виде сложны. Несмотря на то, что по переменным $x$, $y$, $z$ их целевые функции достаточно гладкие, функции, задающие ограничения, таким свойством могут уже не обладать, в частности, это касается функций $g_1$ и $g_2$. То же самое можно сказать и для переменных φ, ψ в задаче (2.5), и для переменных φ, ψ, $r$ в задаче (2.10). По переменным φ, ψ, $r$ задача (2.11) может иметь существенно негладкие и целевую функцию, и функции ограничений.

Разумеется, можно было бы попытаться в любой из задач с помощью определенных процедур провести сглаживание всех используемых в ней функций, после чего, решить сглаженную задачу каким-либо эффективным методом минимизации, например, методом линеаризации или методом последовательного квадратичного программирования. Между тем, сложный вид функций, задающих ограничения (особенно в задачах II.а, II.б (см. ниже)), свидетельствует о нецелесообразности использования существующих процедур сглаживания, ибо существенного упрощения исходной задачи они не дают. Более того, в каждой из задач такой подход может привести к тому, что решение сглаженной задачи может оказаться вне зоны видимости ствола орудия или же довольно далеко от оптимального решения.

Очевидно, что минимизацию целевых функций в задачах (2.5) и (2.10), можно осуществлять как по переменным $(\varphi; \psi)$ и $(\varphi; \psi; r)$, так и по переменным $(x; y)$ и $(x; y; z)$, соответственно. В любом случае найденное решение позволяет определить требуемое оптимальное направление ствола орудия, задаваемое переменными $(\varphi; \psi)$, а в задаче (2.10) – еще и указать ближайшую к цели точку, находящуюся на оптимальной траектории.

Как указывалось выше, в задаче II.б не приходится надеяться на аналитическое определение координат точки $N(x; y; z)$ через переменные $\varphi$ и $\psi$. Поэтому, данную задачу разумно решать по переменным $x$, $y$, $z$. Более того, сейчас мы приведем один из аргументов в пользу того, что и при решении задач (2.5), (2.10) также целесообразно перейти к переменным $x$, $y$ и $z$.

Ввиду указанной выше гладкости по переменным $\varphi$, $\psi$, $r$ целевых функций в задачах (2.5) и (2.10) решать эти задачи можно с помощью какого-либо численного метода первого порядка. Однако, в случае негладкости одной из функций $g_1$ или $g_2$ либо невыпуклости множества $E$, использование данных методов затруднительно, вследствие возможного нарушения условий их сходимости или значительного усложнения решения вспомогательных задач, например, задач проектирования на множество $E$.

В настоящей работе относительно функций $g_1$ и $g_2$, определяющих множество $E$, мы будем лишь предполагать, что они удовлетворяют условию Липшица на некотором прямоугольнике в $\Theta$. Данное ограничение вполне естественно, так как может быть удовлетворено кусочно-линейной аппроксимацией зоны видимости ствола орудия, произведенной, например, на основании результатов фотосъемки.

Учитывая приведенные ограничения на функции $g_1$, $g_2$ и замечания относительно используемых в задачах (2.5), (2.10), (2.11) переменных, решение всех задач будем осуществлять по переменным $x$, $y$, $z$.

Кроме того, ввиду многозначности отображений $h_{\mathrm{I}}^{-1}$ и $h_{\mathrm{II}}^{-1}$ ограничения в каждой из задач целесообразно разбить на два подмножества, соответствующих $h_{k1}^{-1}$ и $h_{k2}^{-1}$, $k \in \{\mathrm{I}; \mathrm{II}\}$ (см. (2.8)). Это разбиение позволяет вместо одной задачи, получить соответствующие ей две подзадачи, которые, очевидно, можно решать независимо друг от друга.

Легко видеть, что каждая из задач (2.5), (2.10), (2.11), (и, следовательно, их подзадач) фактически является задачей проектирования точки $M$ на соответствующее множество, определяемое ограничениями задачи. Это множество для каждой подзадачи, как мы сейчас покажем, может быть описано посредством соответствующего функционального неравенства типа $F_{ij}(x; y; z) \le 0$, $i \in \{\mathrm{I}, \mathrm{II.а}, \mathrm{II.б}\}$, $j \in \{1; 2\}$ ($j$ – номер подзадачи), в котором левая часть удовлетворяет так называемому условию ε-липшицевости. Данное свойство позволит нам в дальнейшем при решении задач использовать методы, разработанные в [2, 3].

Для задач (2.5) и (2.10) функции $F_{\mathrm{I}\,j}$ и $F_{\mathrm{II.a}\,j}$ можно задать равенствами

$$F_{\mathrm{I}\,j}\left(N\right) := \max\left\{\omega_{\mathrm{I}\,j1} g\left(\mathrm{Pr}_1\, h_{\mathrm{I}\,j}^{-1}(N); \mathrm{Pr}_2\, h_{\mathrm{I}\,j}^{-1}(N)\right); \omega_{\mathrm{I}\,j2}\left(\kappa - \mathrm{Pr}_1 N\right)\right\}, \tag{3.1}$$

$$\begin{aligned} F_{\mathrm{II.a}\,j}\left(N\right) := \max\Big\{&\omega_{\mathrm{II.a}\,j1} g\left(\mathrm{Pr}_1\, h_{\mathrm{II},j}^{-1}(N); \mathrm{Pr}_2\, h_{\mathrm{II},j}^{-1}(N)\right); \omega_{\mathrm{II.a}\,j2}\left(\kappa - \mathrm{Pr}_1 N\right); \\ &\omega_{\mathrm{II.a}\,j3}\left(z_{\min} - \mathrm{Pr}_3\, N\right)\Big\}. \end{aligned} \tag{3.2}$$

Здесь и далее $\omega_{ijl}$ – весовые коэффициенты, значения которых зависят от свойств функций, участвующих в задании $F_{i\,j}$. Один способов их выбора будет указан нами позже перед описанием алгоритма решения поставленных задач.

Поскольку в задаче (2.11) ограничения можно формализовать соотношениями

$$\begin{cases} g\left(\mathrm{Pr}_1\, h^{-1}_{\text{П}\,j}(N);\mathrm{Pr}_2\, h^{-1}_{\text{П}\,j}(N)\right)\le 0, H(N)=0, \\ \mathrm{Pr}_1\, N\ge \text{к}, z_{\min}-\mathrm{Pr}_3\, N\le 0, \min\limits_{\lambda\in[0;1]}\{H(\lambda N+(1-\lambda)M)\}\ge 0, \\ \min\limits_{\mu\in[0;1]}\left\{H\left(h_{\text{П}}\left(\mathrm{Pr}_1 h^{-1}_{\text{П}\,j}(N);\mathrm{Pr}_2 h^{-1}_{\text{П}\,j}(N);\mu\,\mathrm{Pr}_3 h^{-1}_{\text{П}\,j}(N)\right)\right)\right\}\ge 0, \end{cases} \tag{3.3}$$

то в качестве функции $F_{\text{П.б}\,j}$ можно взять функцию вида

$$\begin{aligned} F_{\text{П.б}\,j}(N):=\max\Big\{&\omega_{\text{П.б}\,j1}\, g\left(\mathrm{Pr}_1\, h^{-1}_{\text{П}\,j}(N);\mathrm{Pr}_2\, h^{-1}_{\text{П}\,j}(N)\right);\omega_{\text{П.б}\,j2}\left(\text{к}-\mathrm{Pr}_1 N\right); \\ &\omega_{\text{П.б}\,j3}\left|H(N)\right|;-\omega_{\text{П.б}\,j4}\min_{\lambda\in[0;1]}\left\{H\left(\lambda N+(1-\lambda)M\right)\right\}; \\ &-\omega_{\text{П.б}\,j5}\min_{\mu\in[0;1]}\left\{H\left(h_{\text{П}}\left(\mathrm{Pr}_1 h^{-1}_{\text{П}\,j}(N);\mathrm{Pr}_2 h^{-1}_{\text{П},\,j}(N);\mu\,\mathrm{Pr}_3 h^{-1}_{\text{П}\,j}(N)\right)\right)\right\}\Big\}. \end{aligned} \tag{3.4}$$

В системах неравенств (3.3) и функциях $F_{ij}$, $i \in$ {I, П.а, П.б}, $j \in$ {1; 2}, верхние ограничения на значения $x$, $y$ и $z$ в задании множеств $W_i$ (см. (2.4), (2.9)) явно не содержатся, поскольку они обусловлены естественными законами движения снаряда – выраженными уравнениями траекторий полёта снаряда, и учтёнными в условии $N \in h_{\bullet}(E)$. Нижние же границы на значения $x$ и $z$ были наложены дополнительно, исходя из практических соображений, поэтому явно описаны соответствующими компонентами функций максимума в (3.1), (3.2) и (3.4).

Таким образом, каждая из задач $i \in$ {I, П.а, П.б} при каждом $j \in$ {1; 2} приобретает следующий вид

$$\begin{aligned} &\| M-N \| \to \min, \\ &\text{s.t. } F_{i\,j}(N)\le 0. \end{aligned} \tag{3.5}$$

Понятно, что на аналитическое решение этих задач надеяться не приходится. Среди численных методов выбор также небольшой и обусловлен он, главным образом, свойствами функций $F_{i\,j}$.

Оказывается, при довольно естественных предположениях относительно функций $F_{i\,j}$, $i \in$ {I, П.а, П.б}, $j \in$ {1; 2}, каждая из них обладает свойством так называемой ε-липшицевости [4] на соответствующем множестве $W_i$. Приведём соответствующее определение.

Произвольная функция $f : A \to Y$, определенная на подмножестве $A$ нормированного пространства $X$, со значениями в нормированном пространстве $Y$ называется *ε-липшицевой на A*, если для любого $0 < \varepsilon < \varepsilon_0$ найдется такое число $L(\varepsilon) \geq 0$, что при любых $x, y \in A$ справедливо неравенство

$$\|f(x) - f(y)\| \leq L(\varepsilon)\|x - y\| + \varepsilon. \tag{3.6}$$

Множество ε-липшицевых на $A$ функций обозначим символом ε-Lip($A$; $Y$), а нижнюю грань для фиксированной функции $f$ множества функций $L(\varepsilon)$, удовлетворяющих (3.6), (т.е. лучшую оценку величины $L(\varepsilon)$) – символом lip($f$ ; ε) (свойства указанной величины можно найти в [2] и в [4], где для функции $f$ она обозначена, как $l(\varepsilon)$ и $\kappa(\varepsilon)$ соответственно).

Для точных формулировок утверждений об ε-липшицевости функций $F_{ij}$, $i \in$ {I, II.а, II.б}, $j \in$ {1; 2}, а также их обоснований, нам потребуются явные зависимости для $\mathrm{Pr}_1\, h_{k,j}^{-1}(N)$, $\mathrm{Pr}_2\, h_{k,j}^{-1}(N)$, $W_k$, $k \in$ {I; II}, $j \in$ {1; 2}.

Из соотношений (2.3) и (2.7) получаем

$$\varphi(N) := \mathrm{Pr}_1\, h_{k,j}^{-1}(N) = \mathrm{arctg}(y/x), \quad k \in \{\mathrm{I};\mathrm{II}\}, \quad j \in \{1;2\};$$

$$\psi_{\mathrm{I},j}(N) := \mathrm{Pr}_2\, h_{\mathrm{I},j}^{-1}(N) = \frac{\pi}{4} + (-1)^j \left( \frac{\pi}{4} - \frac{1}{2} \arcsin \frac{\sqrt{x^2 + y^2}}{v^2} \right), \quad j \in \{1;2\}; \tag{3.7}$$

$$\begin{aligned} \psi_{\mathrm{II},j}(N) := \mathrm{Pr}_2\, h_{\mathrm{II},j}^{-1}(N) = \mathrm{arctg}\Big( & \left(x^2 + y^2\right)^{-1/2} \times \\ & \times \left[ v^2 + (-1)^j \sqrt{v^4 - \left(x^2 + y^2 + 2v^2 z\right)} \right] \Big), \quad j \in \{1;2\}. \end{aligned} \tag{3.8}$$

Стоит отметить, что при некоторых значениях переменных $x$, $y$ и $z$ выражение, стоящее под вторым квадратным корнем в (3.8), может оказаться отрицательным. Но этот факт лишь свидетельствует о невозможности попадания снаряда в точку ($x$, $y$, $z$) ни при каком значении угла ψ при заданной начальной скорости $v_0$.

ε-липшицевость функций $F_{ij}$, $i \in \{\mathrm{I}, \mathrm{II.а}, \mathrm{II.б}\}$, $j \in \{1; 2\}$, формулируется ниже предложениями 1–3; также приводятся несколько лемм об ε-липшицевости и липшицевости некоторых функций, участвующих в задании функций $F_{ij}$. Всюду предполагается, что $g_1, g_2 \in \mathrm{Lip}\big([\theta_1;\theta_2];\mathbb{R}\big)$, $H \in \mathrm{Lip}(\mathbb{R}^3;\mathbb{R})$, $2\varepsilon\in(0;\pi/2]$, а в получении нижеследующих оценок участвует норма $\|\cdot\|_1$.

**Лемма 1.** *Для функции* $g(\varphi, \psi)$*, определяемой равенством* (2.6)*, справедливо включение* $g \in \mathrm{Lip}(\Theta;\mathbb{R})$*, причем имеет место оценка*

$$\mathrm{lip}(g) \le \max\big\{ \mathrm{lip}(g_1);\mathrm{lip}(g_2);1\big\}. \tag{3.9}$$

Доказательство этой и последующих лемм и предложений приведены в приложении.

**Лемма 2.** *Справедливо включение* $\varphi(x; y) \in \mathrm{Lip}\big(W_k;(-\pi/2;\pi/2)\big)$, $k \in \{\mathrm{I}; \mathrm{II}\}$*, причем имеет место оценка*

$$\mathrm{lip}(\varphi) \le \sqrt{2}\kappa^{-1}. \tag{3.10}$$

**Лемма 3.** *Справедливо включение* $\psi_{\mathrm{I}j}(x; y) \in \varepsilon\text{-}\mathrm{Lip}\big(W_{\mathrm{I}};(0;\pi/2)\big)$, $j \in \{1; 2\}$, *причем имеет место оценка*

$$\mathrm{lip}(\psi_{\mathrm{I}j};\varepsilon) \le \left(\sqrt{2}v^2\sqrt{1-\tau^2(2\varepsilon)}\right)^{-1}, \text{ если } 0<2\varepsilon<\pi/2-1, \tag{3.11}$$

$$\mathrm{lip}(\psi_{\mathrm{I}j};\varepsilon) \le \left(\sqrt{2}v^2\right)^{-1}\big(\pi/2-2\varepsilon\big), \text{ если } \pi/2-1\le 2\varepsilon \le \pi/2,$$

*где* $\tau(2\varepsilon)$ *– корень уравнения* $\big(\pi/2-2\varepsilon-\arcsin\tau\big)\sqrt{1-\tau^2} = 1-\tau$ *на интервале* $[0;1)$.

**Предложение 1.** *Справедливо включение* $F_{\mathrm{I}j}(N) \in \varepsilon\text{-}\mathrm{Lip}\big(W_{\mathrm{I}};\mathbb{R}\big)$, $j \in \{1; 2\}$, *причем имеет место оценка*

$$\mathrm{lip}(F_{\mathrm{I}j},\varepsilon) \le \max\Big\{\omega_{\mathrm{I}j1}\,\mathrm{lip}(g)\cdot\Big(\mathrm{lip}(\varphi)+\mathrm{lip}(\psi_{\mathrm{I}j};\varepsilon(\omega_{\mathrm{I}j1}\,\mathrm{lip}(g))^{-1}\Big);\omega_{\mathrm{I}j2}\Big\}, \tag{3.12}$$

*где* $\mathrm{lip}(g)$, $\mathrm{lip}(\varphi)$ *и* $\mathrm{lip}\Big(\psi_{\mathrm{I}j};\varepsilon(\omega_{\mathrm{I}j1}\,\mathrm{lip}(g))^{-1}\Big)$, $j \in \{1; 2\}$, *определяются по* (3.9), (3.10) *и* (3.11)*, соответственно.*

**Лемма 4.** *Справедливо включение* $\psi_{\text{II}\,j}(x;y)\in\varepsilon\text{-Lip}\big(W_{\text{II}};(-\pi/2;\pi/2)\big)$, $j\in\{1;2\}$, *причем имеет место оценка*

$$\mathrm{lip}\big(\psi_{\text{II}\,j};\varepsilon\big)\le\big(\rho_{\text{II}}(2\varepsilon)^{-1}+\beta\sqrt{2}\big)\kappa^{-2}, \tag{3.13}$$

*где* $\beta=v^2+\sqrt{v^4-\kappa^2-2v^2z_{\min}}$.

**Предложение 2.** *Справедливо включение* $F_{\text{II.а}\,j}(N)\in\varepsilon\text{-Lip}\big(W_{\text{II}};\mathbb{R}\big)$, $j\in\{1;2\}$, *причем имеет место оценка*

$$\mathrm{lip}(F_{\text{II.а}\,j},\varepsilon)\le\max\Big\{\omega_{\text{II.а}\,j1}\,\mathrm{lip}(g)\cdot\Big(\mathrm{lip}(\varphi)+\mathrm{lip}\Big(\psi_{\text{II}\,j};\varepsilon(\omega_{\text{II.а}\,j1}\,\mathrm{lip}(g))^{-1}\Big)\Big);\omega_{\text{II.а}\,j2};\omega_{\text{II.а}\,j3}\Big\}, \tag{3.14}$$

*где* $\mathrm{lip}(g)$, $\mathrm{lip}(\varphi)$ и $\mathrm{lip}\Big(\psi_{\text{II}\,j};\varepsilon(\omega_{\text{II.а}\,j1}\,\mathrm{lip}(g))^{-1}\Big)$, $j\in\{1;2\}$, *определяются по* (3.9), (3.10) *и* (3.13), *соответственно.*

**Предложение 3.** *Справедливо включение* $F_{\text{II.б}\,j}(N)\in\varepsilon\text{-Lip}\big(W_{\text{II}};\mathbb{R}\big)$, $j\in\{1;2\}$, *причем имеет место оценка*

$$\begin{aligned}\mathrm{lip}(F_{\text{II.б}\,j};\varepsilon)\le\max\Big\{&\omega_{\text{II.б}\,j1}\,\mathrm{lip}(g)\cdot\Big(\mathrm{lip}(\varphi)+\mathrm{lip}\Big(\psi_{\text{II}\,j};\varepsilon(\omega_{\text{II.б}\,j1}\,\mathrm{lip}(g))^{-1}\Big)\Big);\,\omega_{\text{II.б}\,j2};\\&\omega_{\text{II.б}\,j3}\,\mathrm{lip}(H);\omega_{\text{II.б}\,j4}\,\mathrm{lip}(H);\omega_{\text{II.б}\,j5}\,\mathrm{lip}(H)\big(\rho_{\text{II}}\,\mathrm{lip}(H)/(8\varepsilon)+1\big)\Big\},\end{aligned} \tag{3.15}$$

*где* $\mathrm{lip}(g)$, $\mathrm{lip}(\varphi)$ *и* $\mathrm{lip}\Big(\psi_{\text{II}\,j};\varepsilon(\omega_{\text{II.б}\,j1}\,\mathrm{lip}(g))^{-1}\Big)$, $j\in\{1;2\}$, *определяются по* (3.9), (3.10) *и* (3.13), *соответственно.*

Согласно нашим предположениям $M\in W_i$, $i\in\{\text{I},\text{II}\}$. Стало быть, $F_{i\,j}(M)\ge 0$, $j\in\{1;2\}$. Эти неравенства и предложения 1–3 показывают, что к задаче проектирования (3.5) можно применить алгоритмы, предложенные в работах [2, 3]. Каждый из алгоритмов строит либо последовательность $Q^m$, сходящуюся к одному из нулей функции $F_{i\,j}$, ближайшему (в смысле выбранной нормы) к точке $M$, причем такую, что $F_{ij}(Q^m)>0$ для всех $m=1, 2, \ldots$, либо за конечное число шагов определяет искомую точку $N_*$, в которой $F_{ij}(N_*)=0$. В настоящей работе мы воспользуемся одним из них. Но прежде чем описать его, обратимся к участвующим при построении функций $F_{ij}$ весовым коэффициентам $\omega_{ijl}$, от выбора которых зависит точность получаемых алгоритмом результатов.

Главной особенностью алгоритма является то, что каждое следующее приближение к искомому решению ищется на основе информации о возможном приращении значений функции $F_{ij}$ относительного текущего приближения. Поскольку функция $F_{ij}$ определяется как максимум от нескольких функций-компонент, то нужно учитывать, как соотносятся между собой приращения этих компонент. Игнорируя же данное сопоставление, мы рискуем тем, что генерируемые алгоритмом приближения не смогут за сколь-нибудь разумное время достигнуть критерия его останова $F_{i\,j}(Q_m) < \varepsilon^*$ (по поводу критерия останова см. также замечания 1–2). Причина этого кроется в том, что точность $\varepsilon^*$ для одной из компонент функции $F_{ij}$ может оказаться слишком завышенной, и, как следствие, число итераций алгоритма для ее достижения – колоссально большим. Как известно [5], верхняя граница числа итераций для алгоритмов нулевого порядка (к которым относится и наш алгоритм), как правило, определяется через приращения значений функции, а значит, зная оценки этих приращений, для компонент функции $F_{ij}$ можно попытаться подобрать такие весовые коэффициенты $\omega_{ijl}$, которые обеспечат адекватное общей точности число итераций. Таким образом, задав некоторую величину параметра общей точности $\varepsilon^*$ (на значения функции $F_{ij}$), можно задавать точности определения отдельных компонент функции $F_{ij}$, используя значения соответствующих весовых коэффициентов: точность $l$-й компоненты будет равна $\varepsilon^*/\omega_{ijl}$.

В настоящей статье положим $\omega_{ij1} = 1$, $i \in \{\text{I, II.а, II.б}\}$, $j \in \{1; 2\}$, т.е. общая точность $\varepsilon^*$ соответствует точности первой компоненты функции $F_{ij}$, для остальных же компонент точности скорректированы соответствующими весовыми коэффициентами $\omega_{ijl}$.

Перейдем теперь к подробному описанию алгоритма.

**Шаг 0.** Задаются величины: начальное значение параметра ε-липшицевости $\varepsilon_0 > 0$; скорость полёта снаряда ($v$ или $v_0$); нижняя граница к значений перемен-ной $x$; начальная точка $N^0 = M \in W_{i,j}$; ограничения $g(\varphi; \psi)$ на зону видимости, задаваемые соотношениями (2.6); параметры алгоритма $\gamma, \lambda \in (0; 1)$. Полагается $Q^0 := N^0$, $k := 0$, $m := 0$.

**Шаг 1.** Вычисляется $F_{i,j}(N^k)$, $i \in \{\text{I, II.а, II.б}\}$, $j \in \{1; 2\}$, в зависимости от типа решаемой задачи по формуле (3.1), (3.2) либо (3.4).

*Если* $F_{i,j}(N^k) < \varepsilon_k(1+\gamma)$, *то* последовательно устанавливаются $N^{k+1} := N^k$, $Q^{m+1} := N^k$, $m := m + 1$, и производится переход к шагу 2; *иначе* – переход к шагу 3.

**Шаг 2.** *Если* $F_{ij}(N^k) \le \varepsilon_k$, *то* принимается $\varepsilon_{k+1} := \lambda F_{i,j}(N^k)$, *иначе* устанавливается $\varepsilon_{k+1} := \lambda\varepsilon_k$. Производится переход к шагу 4.

**Шаг 3.** Определяется $N^{k+1}$ по схеме:

- в зависимости от типа решаемой задачи определяется $\text{lip}(F_{ij}, \varepsilon_k)$, $i \in \{\text{I, II.а, II.б}\}$, $j \in \{1; 2\}$, по одной из формул (3.12), (3.14) или (3.15) с использованием (3.9), (3.10), (3.11), (3.13),
- определяется

$$N^{k+1} = \arg\min\left\{F_{i,j}(X) : X \in \text{fr}\, K_k \cap W_{i,j}\right\}, \quad i \in \{\text{I;II.а;II.б}\}, \quad j \in \{1;2\}, \quad (3.16)$$

где

$$K_k = \left\{X \in R^n : \|X - M\| \le \sum_{s=0}^{k} r_s\right\}, \quad r_k = \frac{F_{i,j}(N^k) - \varepsilon_k}{\sqrt{n}\cdot \text{lip}(F_{i,j}, \varepsilon_k)},$$

$n = 2$ для задачи $i = \text{I}$ и $n = 3$ для задач $i = \text{II.а}$ и $i = \text{II.б}$.

- полагается $\varepsilon_{k+1} := \varepsilon_k$ и осуществляется переход к шагу 4.

**Шаг 4.** *Если* $F_{i,j}(N^{k+1}) = 0$ или $F_{i,j}(Q^m) = 0$, *то* принимается $N^{k+1}$ или $Q^m$, соответственно, в качестве решения задачи (3.5); *иначе* устанавливается $k := k + 1$ и осуществляется переход шагу 1.

**Замечания. 1.** В соответствии с основным утверждением о сходимости алгоритма [2] при бесконечном числе точек $N^k$ он гарантирует сходимость $F_{ij}(Q^m) \underset{m\to\infty}{\longrightarrow} 0$. Поэтому, в данном случае в критерий его останова естественно включить условие вида $F_{ij}(Q^m) < \varepsilon^*$, $\varepsilon^* > 0$. К нему можно также добавить условия, устанавливающие точности для некоторых компонент функции $F_{ij}$, или условие вида $\left\| Q^m - Q^{m+1} \right\| < \varepsilon_Q^*$, $\varepsilon_Q^* > 0$.

**2.** Легко видеть, что основные вычислительные затраты приходятся на третий шаг алгоритма и связаны с решением вспомогательной оптимизационной задачи по нахождению точки $N^{k+1}$. Укажем пару рекомендаций по решению этой задачи. Во-первых, поскольку минимизация функции $F_{ij}$, $j \in \{1; 2\}$, осуществляется либо на окружности (задача $i$ = I), либо на двумерной сфере (задачи $i$ = II.а, II.б), то размерность вспомогательной задачи может быть понижена на единицу переходом к полярным либо к сферическим координатам. Во-вторых, её необязательно решать точно. Поскольку значение $F_{i,j}(N^k)$ на каждом шаге сравнивается с $\varepsilon_k$ и с $\varepsilon_k(1 + \gamma)$, то вместо точки $N^{k+1}$, определяемой по формуле (3.16), можно найти такую точку $\tilde{N}^{k+1}$, что $F_{i,j}(N^{k+1}) \le F_{i,j}(\tilde{N}^{k+1}) < F_{i,j}(N^{k+1}) + \delta_k$, где $\delta_k < \varepsilon_k \gamma$.

**3.** В задаче II.б для вычисления значений функции $F_{\text{II.б}}$ (4-я компонента) необходимо решать задачу одномерной минимизации функции $f(\lambda) = H\left(\lambda N + (1-\lambda) M\right)$. Она может быть приближенно решена каким-либо стандартным методом. Аналогичное замечание касается также и 5-й компоненты функции $F_{\text{II.б}}$.

## 4. ЧИСЛЕННЫЕ ПРИМЕРЫ

Тестирование предложенного в работе метода было проведено на нескольких примерах, в каждом из которых варьировались значения: точности $\varepsilon^*$, которое считалось равным начальному значению $\varepsilon_0$ параметра параметра $\varepsilon$-липшицевости, положение цели $M$, а также вид множества $W_i$, $i \in \{\text{I}, \text{II}\}$. Вспомогательная задача минимизации по вычислению следующей итерационной точки $N^{k+1}$ (шаг 3.2 алгоритма) решалась методом равномерного перебора.

Во всех примерах были зафиксированы следующие значения параметров: $v_0 = 180$ м/с; $\kappa = 100$ м; $\gamma = \lambda = 0{,}5$. Значения варьируемых параметров представлены в таблицах 1–3. При этом использовались обозначения:

- $M_1 = (110; 0; 20)$, $M_2 = (2700; 0; -10)$ – положение цели $M$ (для задачи I третья координата была установлена равной 0);
- $E_1 = \left\{ (\varphi; \psi) : \varphi \in \left[ 0; \pi/2 - 10^{-3} \right], 7\pi/36 \le \psi \le 8\pi/36 \right\}$,

  $E_2 = \left\{ (\varphi; \psi) : \varphi \in \left[ 0; \pi/2 - 10^{-3} \right], (4 + \sin\varphi)\pi/36 \le \psi \le (1 + \sin\varphi)\pi/9 \right\}$;
- $(x_N, y_N, z_N)$ – решение задачи;
- $(\varphi, \psi)$ – соответствующие $(x_N, y_N, z_N)$ углы наклона ствола орудия;
- $k_{\text{общ}}$ – общее число итераций алгоритма;
- $t$ – примерное время работы программы в секундах.

Значение $z_{\min}$ в задачах II.а и II.б принималось равным –10; множество $D$ в задаче II.б задавалось функцией

$$H(x, y, z) = \min \{ \max \{90 - x; x - 130; -10 - y; y - 30; z - 20\}; z + 10\}.$$

Весовые коэффициенты в функциях $F_{\text{I}\,j}$, $F_{\text{II.а}\,j}$ и $F_{\text{II.б}\,j}$ (см. (3.1), (3.2) (3.4)) были установлены следующими:

I) $\omega_{\text{I}\,j1} = 1$, $\omega_{\text{I}\,j2} = 0{,}01$;

II.а) $\omega_{\text{II.а}\,j1} = 1$, $\omega_{\text{II.а}\,j2} = 0{,}01$, $\omega_{\text{II.а}\,j3} = 0{,}01$;

II.б) $\omega_{\text{II.б}\,j1} = 1$, $\omega_{\text{II.б}\,j2} = 0{,}01$, $\omega_{\text{II.б}\,j3} = 0{,}001$, $\omega_{\text{II.б}\,j4} = 0{,}001$, $\omega_{\text{II.б}j5} = 0{,}001$.

Расчет производился на ПК с 4-ядерным ЦПУ IntelCorei3-4020, 1.50 ГГц.

Результаты вычислений приведены в таблицах 1–3. Из них видно, что в задаче I $\varepsilon$-решение находится достаточно быстро. В задаче II.а и особенно в задаче II.б при $\varepsilon_0$ более 0,05 наблюдается существенный рост числа итераций. (Прочерк в таблице 3 означает, что метод не привел к решению в течении 200 сек). Это объясняется тем, что для функций $F_{\text{II.а,б}}$ с уменьшением $\varepsilon$ возрастает оценка $L(\varepsilon)$. Однако, за конечное число шагов (когда оценка $L(\varepsilon)$ разумно велика) все же удается иногда отсечь существенное множество точек, в котором заведомо нет ну-лей функции $F_{ij}$, и таким образом сузить область $W_i$. Разумеется, объем вычислений также можно сократить, по возможности прибегая к распараллеливанию некоторых частей алгоритма. В частности, это можно сделать при решении вспомогательных оптимизационных задач методом равномерного перебора.

## СПИСОК ЛИТЕРАТУРЫ

**Таблица 1. Результаты расчётов для задачи I**

| $M$ | $E$ | ε | $x_N$, м | $y_N$, м | φ, ° | ψ, ° | $k_{общ}$ | $t$, с |
|---|---|---|---|---|---|---|---|---|
| $M_1$ | $E_1$ | 0,1 | 2818,6 | 0,0 | 0,0 | 29,3 | 300 | 0,146 |
| | | 0,05 | 2976,8 | 0,0 | 0,0 | 32,1 | 353 | 0,170 |
| | | 0,01 | 3081,4 | 0,0 | 0,0 | 34,4 | 575 | 0,306 |
| | $E_2$ | 0,1 | 1523,4 | –152,1 | –5,7 | 13,8 | 301 | 0,419 |
| | | 0,05 | 1834,6 | –91,4 | –2,9 | 16,9 | 442 | 0,822 |
| | | 0,01 | 2068,1 | –20,6 | –0,6 | 19,4 | 891 | 3,514 |
| $M_2$ | $E_1$ | 0,1 | 2820,3 | 0,0 | 0,0 | 29,3 | 40 | 0,034 |
| | | 0,05 | 2977,1 | 0,0 | 0,0 | 32,2 | 144 | 0,094 |
| | | 0,01 | 3081,4 | 0,0 | 0,0 | 34,4 | 331 | 0,213 |
| | $E_2$ | 0,1 | 2606,4 | 48,4 | 1,1 | 26,0 | 36 | 0,031 |
| | | 0,05 | 2447,8 | 138,6 | 3,2 | 24,0 | 158 | 0,086 |
| | | 0,01 | 2316,7 | 217,0 | 5,4 | 22,4 | 456 | 0,244 |

**Таблица 2. Результаты расчётов для задачи II.a**

| $M$ | $E$ | ε | $x_N$, м | $y_N$, м | $z_N$, м | φ, ° | ψ, ° | $k_{общ}$ | $t$, с |
|---|---|---|---|---|---|---|---|---|---|
| $M_1$ | $E_1$ | 0,1 | 100,0 | 0,1 | 58,0 | 0,1 | 31,0 | 1889 | 1,249 |
| | | 0,05 | 100,0 | 1,8 | 65,4 | 1,0 | 34,1 | 4182 | 2,244 |
| | | 0,01 | 102,0 | 2,1 | 72,1 | 1,2 | 36,1 | 26297 | 5,227 |
| | $E_2$ | 0,1 | 108,7 | –0,5 | 25,9 | –0,3 | 14,3 | 614 | 1,073 |
| | | 0,05 | 106,8 | –1,0 | 31,4 | –0,6 | 17,3 | 2444 | 1,495 |
| | | 0,01 | 104,8 | –0,3 | 35,9 | –0,2 | 19,8 | 17619 | 8,854 |
| $M_2$ | $E_1$ | 0,1 | 2730,7 | 0,0 | 48,3 | 0,0 | 29,3 | 8478 | 4,064 |
| | | 0,05 | 2771,2 | 0,0 | 121,4 | 0,0 | 32,2 | 33069 | 12,750 |
| | | 0,01 | 2803,2 | –0,1 | 173,5 | –0,0 | 34,4 | 213343 | 103,261 |
| | $E_2$ | 0,1 | 2502,2 | 107,1 | –8,6 | 2,5 | 24,4 | 27281 | 5,952 |
| | | 0,05 | 2336,3 | 826,1 | –4,3 | 19,5 | 24,2 | 244288 | 131,227 |
| | | 0,01 | 2324,5 | 663,2 | –2,8 | 15,9 | 23,4 | 313626 | 163,866 |

**Таблица 3. Результаты расчётов для задачи II.б**

| $M$ | $E$ | ε | $x_N$, м | $y_N$, м | $z_N$, м | φ, ° | ψ, ° | $k_{общ}$ | $t$, с |
|---|---|---|---|---|---|---|---|---|---|
| $M_1$ | $E_1$ | 0,1 | 100,0 | 0,0 | 58,0 | 0,0 | 31,0 | 2748 | 9,753 |
| | | 0,05 | 100,2 | 16,0 | 62,6 | 9,1 | 32,6 | 7308 | 14,240 |
| | $E_2$ | 0,1 | 108,7 | –0,5 | 25,9 | –0,3 | 14,3 | 1072 | 11,303 |
| | | 0,05 | 106,8 | –1,0 | 31,4 | –0,6 | 17,3 | 4258 | 11,179 |
| $M_2$ | $E_1$ | 0,1 | 2730,7 | 0,0 | 48,3 | 0,0 | 29,3 | 12746 | 19,626 |
| | | 0,05 | - | - | - | - | - | - | - |
| | $E_2$ | 0,1 | 2502,2 | 107,1 | –8,6 | 2,5 | 24,4 | 47995 | 88,068 |
| | | 0,05 | 2437,0 | 145,7 | –10,0 | 3,4 | 23,6 | 69230 | 121,723 |

## 5. ПРИЛОЖЕНИЕ

Приведем доказательства лемм и предложений из п.3. При этом, нигде особо не оговаривая, мы будем пользоваться неравенством треугольника и свойствами функций максимума. Также нам понадобится следующая вспомогательная лемма.

**Лемма П.** *Справедливы включения* $\sqrt{x^2+y^2} \in \mathrm{Lip}(W_k;\mathbb{R})$ *и* $\arcsin x \in \varepsilon\text{-}\mathrm{Lip}([0;1];[0;\pi/2])$*, причём по норме* $\|\cdot\|_1$ *имеют место оценки*

$$\mathrm{lip}\left(\sqrt{x^2+y^2}\right) \le \sqrt{2},$$

$$\mathrm{lip}(\arcsin x;\varepsilon) \le \left(1-\tau^2(\varepsilon)\right)^{-1/2}, \text{ если } 0<\varepsilon<\pi/2-1,$$

$$\mathrm{lip}(\arcsin x;\varepsilon) \le \pi/2-\varepsilon, \text{ если } \pi/2-1 \le \varepsilon \le \pi/2,$$

*где* $\tau(\varepsilon)$ *– минимальный корень уравнения* $(\pi/2-\varepsilon-\arcsin\tau)\sqrt{1-\tau^2} = 1-\tau$ *на отрезке* $[0;1]$ *(или,что то же самое,– единственный корень этого же уравнения на интервале* $[0;1)$*).*

Для краткости примем эту лемму без доказательства. ■

**Доказательство леммы 1.** Для компонент функции $g$ вида (2.6) справедливы следующие очевидные оценки:

$$\left|(\theta_1-\varphi_1)-(\theta_1-\varphi_2)\right| = 1\cdot\left|\varphi_1-\varphi_2\right|, \quad \left|(\varphi_1-\theta_2)-(\varphi_2-\theta_2)\right| = 1\cdot\left|\varphi_1-\varphi_2\right|,$$

$$\left|\left(g_1(\varphi_1)-\psi_1\right)-\left(g_1(\varphi_2)-\psi_2\right)\right|\le\left|g_1(\varphi_1)-g_1(\varphi_2)\right|+\left|\psi_1-\psi_2\right|\le \operatorname{lip}(g_1)\left|\varphi_1-\varphi_2\right|+\left|\psi_1-\psi_2\right|=$$
$$=\max\left\{\operatorname{lip}(g_1);1\right\}\cdot\left(\left|\varphi_1-\varphi_2\right|+\left|\psi_1-\psi_2\right|\right),$$

$$\left|\left(\psi_1-g_2(\varphi_1)\right)-\left(\psi_2-g_2(\varphi_2)\right)\right|\le\max\left\{\operatorname{lip}(g_2);1\right\}\cdot\left(\left|\varphi_1-\varphi_2\right|+\left|\psi_1-\psi_2\right|\right).$$

Отсюда для самой функции $g$ получаем

$$\left|g(\varphi_1;\psi_1)-g(\varphi_2;\psi_2)\right|\le\max\left\{\operatorname{lip}(g_1);\operatorname{lip}(g_2);1\right\}\cdot\left(\left|\varphi_1-\varphi_2\right|+\left|\psi_1-\psi_2\right|\right),$$

т.е. $g\in\operatorname{Lip}(\Theta;\mathbb{R})$ и $\operatorname{lip}(g)\le\max\left\{\operatorname{lip}(g_1);\operatorname{lip}(g_2);1\right\}$ по норме $\|\cdot\|_1$. ■

**Доказательство леммы 2.** Используя дифференцируемость функции $\operatorname{arctg}\dfrac{y}{x}$, имеем оценку:

$$\left|\varphi(x_1,y_1)-\varphi(x_2,y_2)\right|\le\left(\left|x_1-x_2\right|+\left|y_1-y_2\right|\right)\sup_{x\in W_k}\frac{\sqrt{2}}{|x|}.$$

Из соотношений (2.4), (2.9), задающих множества $W_k$, $k\in\{\mathrm{I};\mathrm{II}\}$, следует, что $|x|\ge\kappa$, стало быть,

$$\left|\varphi(x_1,y_1)-\varphi(x_2,y_2)\right|\le\sqrt{2}\kappa^{-1}\left(\left|x_1-x_2\right|+\left|y_1-y_2\right|\right).$$

Т.е. имеет место оценка (3.10). ■

**Доказательство леммы 3.** Для обеих ($j\in\{1;2\}$) функций $\psi_{\mathrm{I}j}$ вида (3.7) очевидно, что

$$\left|\psi_{\mathrm{I}j}(x_1,y_1)-\psi_{\mathrm{I}j}(x_2,y_2)\right|\le\frac{1}{2}\left|\arcsin\frac{\sqrt{x_1^2+y_1^2}}{v^2}-\arcsin\frac{\sqrt{x_2^2+y_2^2}}{v^2}\right|.$$

Применяя при $0<\varepsilon_1<\pi/2-1$ лемму П к выражению в правой части этого неравенства, получаем следующую цепочку неравенств:

$$\left|\psi_{\mathrm{I}j}(x_1,y_1)-\psi_{\mathrm{I}j}(x_2,y_2)\right|\le\frac{1}{2}\left(\frac{1}{\sqrt{1-\tau^2(\varepsilon_1)}}\cdot\left|\frac{\sqrt{x_1^2+y_1^2}}{v^2}-\frac{\sqrt{x_2^2+y_2^2}}{v^2}\right|+\varepsilon_1\right)\le$$
$$\le\frac{1}{2}\cdot\frac{1}{\sqrt{1-\tau^2(\varepsilon_1)}}\cdot\frac{\sqrt{2}}{v^2}\left(\left|x_1-x_2\right|+\left|y_1-y_2\right|\right)+\frac{\varepsilon_1}{2}.$$

Полагая $\varepsilon=\varepsilon_1/2$, находим

$$\left|\psi_{\mathrm{I}j}(x_1,y_1)-\psi_{\mathrm{I}j}(x_2,y_2)\right|\le\left(\sqrt{2}v^2\sqrt{1-\tau^2(2\varepsilon)}\right)^{-1}\left(|x_1-x_2|+|y_1-y_2|\right)+\varepsilon\ .$$

Аналогичные рассуждения проводятся и для $\pi/2-1\le\varepsilon_1\le\pi/2$. ■

**Доказательство предложения 1.**

Рассмотрим первую компоненту функции максимума в (3.1). Учитывая липшицевость функции $g$ (см. лемму 1), а также вид компонент отображения $h_{\mathrm{I}j}^{-1}$, имеем

$$\left|g\left(\varphi(N_1);\psi_{\mathrm{I}j}(N_1)\right)-g\left(\varphi(N_2);\psi_{\mathrm{I}j}(N_2)\right)\right|\le\mathrm{lip}(g)\cdot\left(\left|\varphi(N_1)-\varphi(N_2)\right|+\left|\psi_{\mathrm{I}j}(N_1)-\psi_{\mathrm{I}j}(N_2)\right|\right).$$

Далее, используя липшицевость функции $\varphi$ (см. лемму 2) и $\varepsilon$-липшицевость функций $\psi_{\mathrm{I}j}$ (см. лемму 3), соответственно, получаем

$$\left|g\left(\varphi(N_1);\psi_{\mathrm{I}j}(N_1)\right)-g\left(\varphi(N_2);\psi_{\mathrm{I}j}(N_2)\right)\right|\le$$
$$\le\mathrm{lip}(g)\cdot\left(\mathrm{lip}(\varphi)\|N_1-N_2\|_1+\mathrm{lip}(\psi_{\mathrm{I}j},\varepsilon_1)\|N_1-N_2\|_1+\varepsilon_1\right),$$

где $\mathrm{lip}(\varphi)$ и $\mathrm{lip}\left(\psi_{\mathrm{I},j};\varepsilon_1\right)$, $j\in\{1;2\}$, определяются по (3.10) и (3.11), соответственно.

Положив $\varepsilon=\omega_{\mathrm{I}j1}\varepsilon_1\,\mathrm{lip}(g)$, убеждаемся в $\varepsilon$-липшицевости первой компоненты функции максимума в (3.1):

$$\omega_{\mathrm{I}j1}\left|g\left(\varphi(N_1);\psi_{\mathrm{I}j}(N_1)\right)-g\left(\varphi(N_2);\psi_{\mathrm{I}j}(N_2)\right)\right|\le$$
$$\le\omega_{\mathrm{I}j1}\,\mathrm{lip}(g)\cdot\left(\mathrm{lip}(\varphi)+\mathrm{lip}\left(\psi_{\mathrm{I}j},\frac{\varepsilon}{\omega_{\mathrm{I}j1}\,\mathrm{lip}(g)}\right)\right)\|N_1-N_2\|_1+\varepsilon.$$

Вторая компонента функции максимума в (3.1), очевидно, липшицева с постоянной равной $\omega_{\mathrm{I}j2}$.

Таким образом, для функций $F_{\mathrm{I}j}(N)$, $j\in\{1;2\}$ будет иметь место оценка (3.12). ■

**Доказательство леммы 4.** Известно, что $\mathrm{arctg}\,x\in\mathrm{Lip}\left(\mathbb{R};(-\pi/2;\pi/2)\right)$ и $\mathrm{lip}(\mathrm{arctg}\,x)\le1$, потому для каждой ($j\in\{1;2\}$) функции $\psi_{\mathrm{II}j}$ вида (3.8) имеем

$$\left|\psi_{\mathrm{II}\,j}\left(x_1,y_1,z_1\right)-\psi_{\mathrm{II}\,j}\left(x_2,y_2,z_2\right)\right|\le$$

$$\le\left(x_1^2+y_1^2\right)^{-1/2}\cdot\left|\sqrt{v^4-\left(x_1^2+y_1^2+2v^2z_1\right)}-\sqrt{v^4-\left(x_2^2+y_2^2+2v^2z_2\right)}\right|+ \tag{5.1}$$

$$+\left|v^2+\sqrt{v^4-\left(x_1^2+y_1^2+2v^2z_1\right)}\right|\cdot\left|\left(x_1^2+y_1^2\right)^{-1/2}-\left(x_2^2+y_2^2\right)^{-1/2}\right|.$$

Оценим по отдельности выражения, стоящие в правой части этого неравенства.

Из соотношения (2.9), задающего множество $W_{\mathrm{II}}$, получаем

$$\left(x_1^2+y_1^2\right)^{-1/2}\le\kappa^{-1}. \tag{5.2}$$

$$\left|\left(x_1^2+y_1^2\right)^{-1/2}-\left(x_2^2+y_2^2\right)^{-1/2}\right|\le\sqrt{2}\kappa^{-2}\left(\left|x_1-x_2\right|+\left|y_1-y_2\right|\right). \tag{5.3}$$

Поскольку функция $\sqrt{t}$ при $t\ge 0$ является ε-липшицевой с оценкой $\mathrm{lip}(\sqrt{t};\varepsilon)=(4\varepsilon)^{-1}$ (см. [2, 4]), то

$$\left|\sqrt{v^4-\left(x_1^2+y_1^2+2v^2z_1\right)}-\sqrt{v^4-\left(x_2^2+y_2^2+2v^2z_2\right)}\right|\le(4\varepsilon)^{-1}\left|\left(x_1^2-x_2^2\right)+\left(y_1^2-y_2^2\right)+2v^2\left(z_1-z_2\right)\right|+\varepsilon.$$

Отсюда элементарно приходим к неравенству

$$\left|\sqrt{v^4-\left(x_1^2+y_1^2+2v^2z_1\right)}-\sqrt{v^4-\left(x_2^2+y_2^2+2v^2z_2\right)}\right|\le$$

$$\le(4\varepsilon)^{-1}\left(\left|x_1-x_2\right|\left(\left|x_1\right|+\left|x_2\right|\right)+\left|y_1-y_2\right|\left(\left|y_1\right|+\left|y_2\right|\right)+2v^2\left|z_1-z_2\right|\right)+\varepsilon.$$

Снова учитывая (2.9), получаем соотношения

$$\left|x\right|,\left|y\right|\le\rho_{\mathrm{II}},\ v^2\le\rho_{\mathrm{II}},\ 0<\kappa\le x,$$

с помощью которых находим

$$\left|\sqrt{v^4-\left(x_1^2+y_1^2+2v^2z_1\right)}-\sqrt{v^4-\left(x_2^2+y_2^2+2v^2z_2\right)}\right|\le\frac{\rho_{\mathrm{II}}}{2\varepsilon}\left(\left|x_1-x_2\right|+\left|y_1-y_2\right|+\left|z_1-z_2\right|\right)+\varepsilon;$$

$$\left|v^2+\sqrt{v^4-\left(x_1^2+y_1^2+2v^2z_1\right)}\right|\le\beta=v^2+\sqrt{v^4-\kappa^2-2v^2z_{\min}}.$$

Подставляя найденные оценки (5.2)–(5.3) в (5.1) и полагая затем $\varepsilon=\varepsilon_1/\kappa$, окончательно получаем

$$\left|\psi_{\mathrm{II}\,j}\left(x_1,y_1,z_1\right)-\psi_{\mathrm{II}\,j}\left(x_2,y_2,z_2\right)\right|\le\kappa^{-2}\left(\rho_{\mathrm{II}}(2\varepsilon)^{-1}+\beta\sqrt{2}\right)\left\|N_1-N_2\right\|_1+\varepsilon.$$

■

**Доказательство предложения 2** проводится полностью аналогично доказательству предложения 1.

**Доказательство предложения 3.** $\varepsilon$-липшицевость первой компоненты и липшицевость второй компоненты вытекает из предложений 1–2. Липшицевость третьей и четвертой компонент следует напрямую из предположения о липшицевости функции $H$. Докажем $\varepsilon$-липшицевость пятой компоненты.

Пользуясь известными оценками функции минимума [6] и указанной липшицевостью функции $H$, находим:

$$\left|\min_{\mu\in[0;1]} H\left(h_{\text{II}}\left(\varphi(N_1);\psi_{\text{II},j}(N_1);\mu r(N_1)\right)\right) - \min_{\mu\in[0;1]} H\left(h_{\text{II}}\left(\varphi(N_2);\psi_{\text{II},j}(N_2);\mu r(N_2)\right)\right)\right| \le$$
$$\le \operatorname{lip}(H)\sum_{i=1}^{3}\max_{\mu\in[0;1]}\left|e_i(N_1)-e_i(N_2)\right|, \tag{5.4}$$

где $N_1(x_1, y_1, z_1)$ и $N_2(x_2, y_2, z_2)$ и $e_i(N)=\Pr_i h_{\text{II}}\left(\varphi(N);\psi_{\text{II},j}(N);\mu r(N)\right)$, $i \in \{1, 2, 3\}$.

Оценим каждое из слагаемых правой части неравенства (5.4).

Пользуясь соотношениями (2.7) и тем, что $\mu \in [0; 1]$, получаем

$$\max_{\mu\in[0;1]}\left|e_1(N_1)-e_1(N_2)\right|=\left|x_1-x_2\right|,\ \max_{\mu\in[0;1]}\left|e_1(N_2)-e_2(N_2)\right|=\left|y_1-y_2\right|, \tag{5.5}$$

$$\max_{\mu\in[0;1]}\left|e_3(N_1)-e_3(N_2)\right|\le\left|r(N_1)\operatorname{tg}\psi_{\text{II},j}(N_1)-r(N_2)\operatorname{tg}\psi_{\text{II},j}(N_2)\right|\max_{\mu\in[0;1]}\mu(1-\mu)+\left|z_1-z_2\right|.$$

Для оценивания первого слагаемого в правой части последнего неравенства воспользуемся доказательством леммы 4 и тем, что $\max\limits_{\mu\in[0;1]}\mu(1-\mu)=1/4$. Имеем

$$\max_{\mu\in[0;1]}\left|e_3(N_1)-e_3(N_2)\right|\le\left(\rho_{\text{II}}(8\varepsilon)^{-1}+1\right)\left\|N_1-N_2\right\|_1+\varepsilon. \tag{5.6}$$

Подставляя оценки (5.5) и (5.6) в (5.4) и полагая затем $\varepsilon = \varepsilon_1 \operatorname{lip}(H)$, получаем

$$\left|\min_{\mu\in[0;1]} H\left(h_{\text{II}}\left(\varphi(N_1);\psi_{\text{II},j}(N_1);\mu r(N_1)\right)\right) - \min_{\mu\in[0;1]} H\left(h_{\text{II}}\left(\varphi(N_2);\psi_{\text{II},j}(N_2);\mu r(N_2)\right)\right)\right| \le$$
$$\le \operatorname{lip}(H)\left(\rho_{\text{II}}\operatorname{lip}(H)(8\varepsilon)^{-1}+1\right)\left\|N_1-N_2\right\|_1+\varepsilon.$$

Объединяя оценки компонент функции (3.4), окончательно приходим к (3.15). ■

# MODELS AND METHODS FOR THREE EXTERNAL BALLISTICS INVERSE PROBLEMS

N.K. ARUTYUNOVA, A.M. DULLIEV, V.I. ZABOTIN

**Abstract**. We consider three problems of selecting optimal gun barrel direction (or those of selecting optimal semi-axis position) when firing an unguided artillery projectile on the assumption that the gun barrel semi-axis can move in a connected nonconvex cone having a non-smooth lateral surface and modelling visibility zone restrictions. In the first problem, the target is in the true horizon plane of the gun, the second and the third problems deal with some region of 3D space. A distinctive feature of the models is that the objective functions are ε-Lipschitz ones. We have constructed a unified numerical method to solve these problems based on an algorithm of projecting a point onto ε-Lipschitz level function set. A computer program has been based on it. series of numerical experiments on each problem has been carried out.

## Introduction

Let us consider a problem of firing an unguided projectile. It is required to minimize the distance from the projectile drop point to the target. This problem is a problem of the external ballistic theory. It has been studied fairly well (see [1]) if the following conditions hold: air resistance is not taken into account, movement of the gun barrel is restricted by the true horizon plane of the gun, and the Earth surface is spherical.

However, in reality, a gun barrel direction can be arbitrarily selected, as a rule, only within some connected nonconvex cone that has non-smooth lateral surface and arises in conditions narrowing the selecting gun barrel direction due to some obstacles.

To restrictions on the selecting gun barrel direction, it is often necessary to add conditions such that the target can lie outside of the true horizon plane of the gun or, possibly, on the surface defining the relief of a landscape. In the latter case, the problem of selecting optimal gun barrel direction becomes much more complicated. The foremost reason for this difficulty is that the minimal Euclidean distance from the target to the projectile drop point does not always correspond to the optimal shot or even close to the optimal one. For example, this lack of correspondence can be observed in case the target and the projectile drop point are separated from each other by some substantial obstacle.

A mathematical model can be described as follows. Suppose $\mathrm{dist}(,)$ is the Euclidean distance, $O$ is the point in which located the gun, $l$ is a ray with origin $O$, $N$ is the projectile drop point, $N = N(l)$, $M$ is the target, $K$ is a cone with vertex $O$, $D'$ is a set in space $\mathbb{R}^3$, and $M \in D'$. It is required

$$\min_{l} \mathrm{dist}(N(l), M),$$
$$\text{subject to} \quad l \in K, \ N(l) \in D'.$$

Naturally, there may be other possible mathematical models with other objective functions. In particular, below we consider a problem in which the objective function is the distance between the projectile trajectory and the point $M$.

We stress that, in this paper, air resistance is not taken into account. At the same time, first, our models that are described and investigated below can be applied for firing an unguided heavy projectile. Secondly, the solutions obtained in the models can be served as base for further finding for a more precise solution and can be refined by appropriate methods.

Below, we formulate and investigate three problems for which a unified algorithm for their solving is constructed. Also, we present numerical results for several test examples; these results were obtained with a computer program implementing the algorithm.

## 1. Formulation of The Problems

Throughout, the gun is modelled by a point $T$, the direction of the gun barrel is modelled by a ray starting from this point, the oblateness of the Earth is not taken account.

Let us choose a Cartesian coordinate system $Oxyz$ with the origin $O$ coinciding with $T$. Suppose the plane $Oxy$ is the horizontal plane at $O$, the axis $Oz$ is directed vertically up, and the gravitational acceleration $\boldsymbol{g}$, vertically down. The directions of the axes $Ox$ and $Oy$ will be given below. Throughout, we use the following notation: $\|\cdot\|_p$ is the $p$-norm on $\mathbb{R}^n$, $n \geqslant 1$ (the subscript $p$ will be omitted for $p = 2$); int and fr are the operators of taking interior and boundary in $\mathbb{R}^n$, respectively; $\mathrm{Pr}_i A$ is the projection of a set $A$ on the $i$-th axis (here, the axes $x$, $y$, and $z$ are denoted by 1, 2, and 3, respectively); given $X \subset \mathbb{R}^n$ and $Y \subset \mathbb{R}$, we denote by $\mathrm{Lip}(X; Y)$ the set of all Lipschitz continuous on $X$ functions $f: X \to Y$ with Lipschitz constant $lip(f)$.

We will assume that the direction (i.e. ray position) of the gun barrel can be freely selected within a closed cone $K$ with the apex $O$ and the cone $K$ contains no vertical rays starting from $O$. Neglecting the length of the gun barrel, suppose the initial projectile velocity vector $\boldsymbol{v}$ has the initial point at $O$ and has a constant length $v_0$. In other words, $\boldsymbol{v} \in S \cap K$, where $S = \{\boldsymbol{v} \mid \|\boldsymbol{v}\| = v_0\} \subset \mathbb{R}^3$ is the sphere of radius $v_0$ centered at $O$. By $v$ denote the fraction $v_0/\sqrt{\|\boldsymbol{g}\|}$.

Throughout, we assume that the selecting direction of the gun barrel is equivalent to the selecting the vector $\boldsymbol{v} = (v_x, v_y, v_z)$ with a given fixed $v_0$. This vector will be described by means of a spherical coordinate system in the form

$$\begin{gathered} v_x = v_0 \cos\psi \cos\varphi, \quad v_y = v_0 \cos\psi \sin\varphi, \quad v_z = v_0 \sin\psi, \\ (\varphi, \psi) \in E \subset \Theta := \{(\varphi, \psi) \mid \varphi \in (-\pi, \pi], \psi \in (-\pi/2, \pi/2)\}. \end{gathered} \tag{1}$$

where the set $E$ is closed and corresponds to the set $S \cap K$, which is given by coordinates $(\varphi, \psi)$.

We assume that the air resistance is negligible. It follows that the trajectory travelled by a projectile is a parabolic trajectory on the gun plane.

In this paper we will consider three problems of selecting optimal gun barrel direction when firing an unguided projectile to the given target. We will suppose that the target $M$ is the point with $y$-coordinate equal to zero and $M$ belongs to the convex set bounded by the paraboloid reached by the projectile (see [1]), and this paraboloid is obtained by union of all the projectile trajectories on the condition of the constraints absence on gun barrel directions. (The case, where the target can lie outside of the paraboloid, and the problem of joint selecting optimal gun barrel direction $(\varphi, \psi)$ and magnitude $v_0$ of the initial projectile velocity are not considered in this paper.)

Among the points lying on the projectile trajectories, we will consider only the points with $x$-coordinate such that $x \geqslant \kappa$, where $\kappa = const$ and $\kappa > 0$. In practical terms, the constant $\kappa$, for example, may correspond to the blast radius of a projectile.

Other constraints on $x$, $y$, and $z$-coordinates are different and depend on the problems that will be considered below.

**Problem I.** Let the point $M$ be the target, lie on the plane $xOy$, and have coordinates $(a, 0)$, $a > \kappa$. Suppose the projectile drop point $N$ also lies on the plane $xOy$ and has coordinates $(x, y)$. For all the trajectories, there are no barriers determined by $\boldsymbol{v} \in S \cap K$. Out of all those trajectories, it is required to choose one for which the Euclidean distance between the points $M$ and $N$ is minimal.

Following [1], the distance $r$ from $O$ to $N(x, y)$ is calculated by

$$r = \sqrt{x^2 + y^2} = v^2 \sin 2\psi. \tag{2}$$

Using (1) and (2), we see that the coordinates $x$ and $y$ are determined from

$$x = r\cos\varphi = v^2 \sin 2\psi \cos\varphi, \quad y = r\sin\varphi = v^2 \sin 2\psi \sin\varphi. \tag{3}$$

Determined by (3) the map $(\varphi,\psi)\mapsto(x,y)$ from $\Theta$ to $\mathbb{R}^2$ is denoted by $h_{\mathrm{I}}$.

Since $v$ is fixed, we obtain the upper bound $\rho_{\mathrm{I}}=v^2$ for the distance travelled by a projectile. Let us introduce notation

$$W_{\mathrm{I}}=\left\{(x,y)\middle|\sqrt{x^2+y^2}\leqslant\rho_{\mathrm{I}},x\geqslant\kappa\right\}. \tag{4}$$

Obviously, the feasible set of the projectile drop points is given by the set $h_{\mathrm{I}}(E)$.

Thus, the gun barrel direction determined by pair $(\varphi,\psi)$ and the corresponding optimal projectile trajectory can be found by solving the problem

$$\begin{aligned}&\| M-N \|^2\to\min,\\&\mathrm{s.t.}\, N\in h_{\mathrm{I}}(E)\cap W_{\mathrm{I}}.\end{aligned} \tag{5}$$

In this paper the set $E$ will be described by means of functional inequalities.

Let $[\theta_1,\theta_2]\subset(-\pi/2,\pi/2)$ and the functions $g_1(\varphi)$ and $g_2(\varphi)$ satisfy the condition

$$-\pi/2<g_1(\varphi)\leqslant g_2(\varphi)<\pi/2,\ \varphi\in[\theta_1,\theta_2].$$

We put

$$E=\{(\varphi,\psi)\in\Theta|\varphi\in[\theta_1,\theta_2],g_1(\varphi)\leqslant\psi\leqslant g_2(\varphi)\},$$

or, in other words,

$$E=\{(\varphi,\psi)\in\Theta|g(\varphi,\psi)\leqslant 0\},$$

where

$$g(\varphi,\psi):=\max\{\theta_1-\varphi,\varphi-\theta_2,g_1(\varphi)-\psi,\psi-g_2(\varphi)\}. \tag{6}$$

Note that Problem I admits a generalization to the case when the target is the set $\mathcal{M}=\{M_1,M_2,\dots,M_n\}$ of the points $M_i=(a_i,b_i)$, $i=1,\dots,n$, $n>1$, on the plane $xOy$. Indeed, the Chebyshev center of $\mathcal{M}$ can be choosen as the target at a single shot, for example.

**Problem II.** Let the point $M(a,0,c)\in\mathbb{R}^3$ be the target. Suppose the target and all the possible projectile drop points lie not below the plane $z=z_{\min}$ ($z_{\min}<0$). As in Problem I, assume that $a>\kappa$. We will distinguish two subproblems, which depend on constraints imposed on the point $M$. These subproblems will be named Problem II.a and Problem II.b, respectively.

**II.a.** Let the point $M$ be an arbitrary point (i.e. $M$ can be in the air) and the projectile trajectory selection is limited only the set $E$. Among all admissible trajectories, it is required to choose the one for which the distance from it to the target $M$ is minimal in comparison with the distance from any remaining admissible trajectory to the target $M$.

**II.b.** Let the point $M$ belong to the boundary $\mathrm{fr}D$ of a set $D$ determined by inequality $H(x,y,z)\leqslant 0$, where $H$ is a continuous function over $\mathbb{R}^3$ and $H(0,0,0)\geqslant 0$. This function defines the relief of a landscape. (It means that the target can be located on a surface of ground/water or some fixed ground object, and the gun, on a surface of ground/water or in the air). We assume that the point $N(x,y,z)$ belongs to $\mathrm{fr}D$, the section of the trajectory from $O$ to $N$ has no common points with $\mathrm{int}D$, and $z\leqslant z_{\min}$ implies $H(x,y,z)\leqslant 0$. Among all admissible trajectories, it is required to choose the one for which the segment $MN$ does not intersect $\mathrm{int}D$ and has minimal length or it is required to be certain that there is no such one. This requirement means that the projectile drop point must be within the line-of-sight of the target and be located maximally close to the target.

Note that the target can lie above ($c>0$) or below ($c<0$) relative to the plane $xOy$ and the optimal value of angle $\psi$ can be less than zero when $c<0$.

Before giving analytical formulations of Problem II.a and Problem II.b, let us make some important remarks.

Unlike Problem I, for Problem II.a and Problem II.b, it is either impossible or very difficult to specify a single-valued map from $\Theta$ to $\mathbb{R}^3$ that analogous to $h_{\mathrm{I}}$. A reason for this difficulty is that the target and the coverage zone in Problems II.a and II.b are located in $\mathbb{R}^3$. Speaking in more detail, in Problem II.a, to each point $(\varphi,\psi)\in E$ there corresponds not one point but a subset of the set of projectile trajectory points; in Problems II.b an analytical determination of coordinates of the point $N(x,y,z)$ for given $(\varphi,\psi)\in E$ may be realized in only special cases that is defined by the set $D$.

It is convenient to consider a single-valued map instead of a set-valued map from $\Theta$ to $\mathbb{R}^3$. This map, denoted below by $h_{\mathrm{II}}\colon \Theta\times\mathbb{R}^+\to\mathbb{R}^3$, takes each $(\varphi,\psi,r)\in\Theta\times\mathbb{R}^+$ to the point $(x,y,z)\in\mathbb{R}^3$ assuming that $(x,y,z)$ lies on the projectile trajectory determined by $(\varphi,\psi)$ and the distance from $O$ to the projection of the point $(x,y,z)$ on the plane $xOy$ is equal to $r$.

According to the theory of external ballistic, the map $h_{\mathrm{II}}$ is determined by

$$h_{\mathrm{II}}(\varphi,\psi,r)=\{(x,y,z)|x=r\cos\varphi, y=r\sin\varphi, z=r\tan\psi-(1+\tan^2\psi)(2v^2)^{-1}r^2\}. \quad (7)$$

Despite single-valuedness of the maps $h_{\mathrm{I}}$ and $h_{\mathrm{II}}$, their inverse maps $h_{\mathrm{I}}^{-1}$ and $h_{\mathrm{II}}^{-1}$ are multivalued: the drop point can be reached by selecting one of the two appropriate values of the angle $\psi\in(-\pi/2,\pi/2)$ for the direction of the gun barrel. At the same time, the value of the angle $\varphi$ is uniquely determined. In constructing the functional constraints for Problem I and II.a,b, the maps $h_{\mathrm{I}}^{-1}$ and $h_{\mathrm{II}}^{-1}$ play an important role. A reason, on which this constructing is possible, is that the maps $h_{\mathrm{I}}^{-1}$ and $h_{\mathrm{II}}^{-1}$ can be parameterized by single-valued maps. Indeed, it follows from (3), (7) that there exist the single-valued maps $h_{\mathrm{I}j}^{-1}$, $h_{\mathrm{II}j}^{-1}$, $j\in\{1,2\}$ for which

$$\begin{aligned}
&h_k^{-1}=h_{k1}^{-1}\cup h_{k2}^{-1},\ \mathrm{dom}(h_{kj}^{-1})=\mathrm{dom}\,(h_k^{-1}),\ k\in\{\mathrm{I},\mathrm{II}\},\ j\in\{1,2\};\\
&\mathrm{Pr}_2\mathrm{Im}(h_{\mathrm{I}1}^{-1})=(0,\pi/4],\ \mathrm{Pr}_2\mathrm{Im}(h_{\mathrm{I}2}^{-1})=[\pi/4,\pi/2),\\
&\mathrm{Pr}_2\mathrm{Im}(h_{\mathrm{II}1}^{-1})=(-\pi/2,\pi/2),\ \mathrm{Pr}_2\mathrm{Im}(h_{\mathrm{II}2}^{-1})=(0,\pi/2).
\end{aligned} \quad (8)$$

The explicit expressions for $h_{kj}^{-1}$ will be given below.

Since all the projectile drop points lie not below the plane $z=z_{\min}$ and using (7), the value $r$ can be bounded from above by

$$\rho_{\mathrm{II}}=v\sqrt{v^2-2z_{\min}}.$$

Taking into account this bound and the envelope of the family of the possible projectile trajectories, we define the set

$$W_{\mathrm{II}}=\left\{(x,y,z)|\sqrt{x^2+y^2}\leqslant\rho_{\mathrm{II}},\ x\geqslant\kappa,\ z_{\min}\leqslant z\leqslant\tfrac{1}{2}\left(v^2-\tfrac{x^2+y^2}{v^2}\right)\right\}, \quad (9)$$

which is the set of all the points reached by a projectile.

In view of the above remarks, we can now give an analytical formulation of Problem II.a in the form

$$\begin{aligned}
&\|M-N\|^2\to\min,\\
&\mathrm{s.t.}\,N\in h_{\mathrm{II}}(E)\cap W_{\mathrm{II}}.
\end{aligned} \quad (10)$$

For Problem II.b we have such an analytical formulation

$$\begin{aligned}
&\|M-N\|^2\to\min,\\
&\mathrm{s.t.}\,N\in h_{\mathrm{II}}(E)\cap W_{\mathrm{II}}\cap\mathrm{fr}D,\ \{\lambda N+(1-\lambda)M|\lambda\in[0,1]\}\cap\mathrm{int}D=\varnothing,\\
&\exists j\in\{1,2\}\colon h_{\mathrm{II}}\big(\mathrm{Pr}_1h_{\mathrm{II}j}^{-1}(N),\mathrm{Pr}_2h_{\mathrm{II}j}^{-1}(N),[0,\mathrm{Pr}_3h_{\mathrm{II}j}^{-1}(N)]\big)\cap\mathrm{int}D=\varnothing.
\end{aligned} \quad (11)$$

Obviously, the last two equations in (11) mean that the section of the trajectory from $O$ to $N$ and the segment $MN$ does not intersect $\mathrm{int}D$, respectively.

By $N_*$ denote a point that belongs to the solution set of problem (11).

Unfortunately, problem (11) admits the existence of common points, in addition to $N_*$, between the segment $MN$ or the section of the trajectory from $O$ to $N_*$ and $\mathrm{fr}D$. Of course, adding conditions to the constrains of problem (11), we can eliminate this situation. However, as it is easy to show, a formalization of these conditions does not allow us to correctly use any optimization method because the feasible set is not closed, moreover, an attempt to take its closure again leads to problem (11). To solve this problem we propose the following procedure. It is necessary to solve problem (11) and if its solution $N_*$ does not satisfy the condition given at the beginning of the paragraph, we have to find another solution of problem (11) in the union of sufficiently small closed "ring like" neighborhoods of the points $(\mathrm{Pr}_1 h_{\mathrm{II}1}^{-1}(N_*), \mathrm{Pr}_2 h_{\mathrm{II}1}^{-1}(N_*))$ and $(\mathrm{Pr}_1 h_{\mathrm{II}2}^{-1}(N_*), \mathrm{Pr}_2 h_{\mathrm{II}2}^{-1}(N_*))$, respectively.

## 2. Minimizing The Objective Functions

Problems (5), (10), and (11) are complex in general case. Despite the fact that the objective functions for these problems are smooth at every point $(x, y, z)$ in their domains, the functions describing the feasible set can be non-differentiable. In particular, this applies to the functions $g_1$ and $g_2$. Similarly, we can assert that these statements hold for problems (5) and (10) with respect to the variables $\varphi$, $\psi$ and the variables $\varphi$, $\psi$, $r$, respectively. Moreover, problem (11) with respect to the variables $\varphi$, $\psi$, $r$ may have both objective function and feasible set functions that are non-differentiable.

Surely, in each problem (5), (10), and (11), we could try smoothing all its functions and then solve the corresponding smoothed problem by invoking some optimization method, for example, linearization method or sequential quadratic programming method. At the same time, a complicated form of the feasible set functions, especially in Problems II.a, II.b, (see below), indicates that using smoothing procedures is not useful, because the ones do not give substantial simplification of the original problem. Moreover, after this smoothing, the solution of any smoothed problem may be outside the visibility zone $E$ for the gun barrel or may be far from the optimal solution.

It is obvious that we can minimize the objective function in problem (5) either with respect to the variables $(\varphi, \psi)$, or with respect to the variables $(x, y)$. Similarly, this holds for problem (10) with respect to the variables $(\varphi, \psi, r)$ and $(x, y, z)$, respectively. In any case, the found optimal solution allows us to find right away the required optimal gun barrel direction, which is determined by $(\varphi, \psi)$, and, for problem (10), to find, in addition, a point closest to the target and lying on the optimal trajectory.

As above, in Problem II.b, we should not rely on an analytical determination of coordinates of the point $N(x, y, z)$ by using the variables $\varphi$, $\psi$. Therefore, it is reasonable to solve this problem with respect to the variables $x$, $y$, $z$. Moreover, we give an argument, which confirms that problems (5) and (10) should be also solved with respect to the variables $x$, $y$, and $z$.

Thanks to smoothness of the objective functions for problems (5) and (10) with respect to the variables $\varphi$, $\psi$, $r$, we can solve these problems using a first-order optimization algorithm. However if the functions $g_1$ or $g_2$ are non-differentiable or the set $E$ is non-convex, then the use of this algorithm is difficult, because the conditions for its convergence may be not satisfied or because solving the auxiliary problems, for example, as finding a projection onto $E$, is also difficult.

In this paper we will suppose that the functions $g_1$ and $g_2$ are only Lipschitz continuous on some solid rectangle contained in $\Theta$. This condition is completely natural, since it can be satisfied, if the visibility zone $E$ has been piecewise linearly approximated on the basis of the photography results.

Following the conditions for the functions $g_1$ and $g_2$, and following the above remarks about the variables in problems (5), (10), and (11), we will solve all problems with respect to $x$, $y$, $z$.

In addition, since the maps $h_{\mathrm{I}}^{-1}$ and $h_{\mathrm{II}}^{-1}$ are multivalued, it follows that the feasible sets in each problems should be presented as the union two subsets that correspond to $h_{k1}^{-1}$ and $h_{k2}^{-1}$, $k \in$ $\{\mathrm{I}, \mathrm{II}\}$ (see (8)). This implies that instead of one problem either (5), or (10), or (11), we have two corresponding subproblems, which may be solved independently of each other.

Clearly, each problem (5), (10), or (11) (and, therefore, its subproblems) is a problem of finding a projection the point $M$ onto the feasible set. We claim that the feasible set for each problem can be defined by the corresponding functional inequality in the form

$$F_{ij}(x, y, z) \leqslant 0, \ i \in \{\mathrm{I}, \mathrm{II.a}, \mathrm{II.b}\}, \ j \in \{1,2\},$$

where $j$ is the subproblem number and the left-hand side of inequality satisfies so-called the $\varepsilon$-Lipschitz condition. This condition will allow us to use the methods proposed in [2, 3].

For problems (5) and (10), the functions $F_{\mathrm{I}j}$ and $F_{\mathrm{II.a}\,j}$ can be defined by

$$F_{\mathrm{I}j}(N) := \max\left\{\omega_{\mathrm{I}j1}\, g\left(\mathrm{Pr}_1 h_{\mathrm{I}j}^{-1}(N), \mathrm{Pr}_2 h_{\mathrm{I}j}^{-1}(N)\right), \omega_{\mathrm{I}j2}(\kappa - \mathrm{Pr}_1 N)\right\}, \tag{12}$$

$$\begin{aligned} F_{\mathrm{II.a}\,j}(N) := \max\{ \quad & \omega_{\mathrm{II.a}\,j1}\, g\left(\mathrm{Pr}_1 h_{\mathrm{II}j}^{-1}(N), \mathrm{Pr}_2 h_{\mathrm{II}j}^{-1}(N)\right), \\ & \omega_{\mathrm{II.a}\,j2}(\kappa - \mathrm{Pr}_1 N), \omega_{\mathrm{II.a}\,j3}(z_{\min} - \mathrm{Pr}_3 N)\}. \end{aligned} \tag{13}$$

Here and in what follows, by $\omega_{ijl}$ we denote the weighting coefficients whose values depend on the properties of the functions defining $F_{ij}$. Later on, before presenting the algorithm for solving Problems I, II.a, II.b, we will propose a technique for choosing $\omega_{ijl}$.

Since the constraints for problem (11) can be formalized by

$$\begin{cases} g\left(\mathrm{Pr}_1 h_{\mathrm{II}j}^{-1}(N), \mathrm{Pr}_2 h_{\mathrm{II}j}^{-1}(N)\right) \leqslant 0, H(N) = 0, \mathrm{Pr}_1 N \geqslant \kappa, \min\limits_{\lambda \in [0,1]} \{H(\lambda N + (1-\lambda)M)\} \geqslant 0, \\ \min\limits_{\mu \in [0,1]} \left\{H\left(h_{\mathrm{II}}\left(\mathrm{Pr}_1 h_{\mathrm{II}j}^{-1}(N), \mathrm{Pr}_2 h_{\mathrm{II}j}^{-1}(N), \mu \mathrm{Pr}_3 h_{\mathrm{II}j}^{-1}(N)\right)\right)\right\} \geqslant 0, z_{\min} - \mathrm{Pr}_3 N \leqslant 0, \end{cases} \tag{14}$$

we can define the function $F_{\mathrm{II.b}\,j}$ by the formula

$$\begin{aligned} F_{\mathrm{II.b}\,j}(N) := \max\{ \quad & \omega_{\mathrm{II.b}\,j1}\, g\left(\mathrm{Pr}_1 h_{\mathrm{II}j}^{-1}(N), \mathrm{Pr}_2 h_{\mathrm{II}j}^{-1}(N)\right), \omega_{\mathrm{II.b}\,j2}(\kappa - \mathrm{Pr}_1 N), \\ & \omega_{\mathrm{II.b}\,j3}|H(N)|, -\omega_{\mathrm{II.b}\,j4} \min_{\lambda \in [0,1]} \{H(\lambda N + (1-\lambda)M)\}, \\ & -\omega_{\mathrm{II.b}\,j5} \min_{\mu \in [0,1]} \left\{H\left(h_{\mathrm{II}}\left(\mathrm{Pr}_1 h_{\mathrm{II}j}^{-1}(N), \mathrm{P}Pr_2 h_{\mathrm{II}j}^{-1}(N), \mu \mathrm{Pr}_3 h_{\mathrm{II}j}^{-1}(N)\right)\right)\right\}\Bigg\}. \end{aligned} \tag{15}$$

The upper bounds for the variables $x$, $y$ and $z$ in the definition of the sets $W_i$ (see (4), (9)) are implicitly contained in (12)–(15), since they follow from projectile motion equations and from the condition $N \in h_{\cdot\cdot}(E)$. But the lower bounds for the variables $x$ and $z$ are explicitly contained in (12)–(15), since they were introduced additionally, based on practical considerations.

Thus, each Problem $i$, $i \in \{\mathrm{I}, \mathrm{II.a}, \mathrm{II.b}\}$, for each $j \in \{1,2\}$ has the form

$$\begin{aligned} & \| M - N \|^2 \to \min, \\ & \text{s.t. } F_{ij}(N) \leqslant 0. \end{aligned} \tag{16}$$

Obviously, we should not rely on an analytical methods for solving any problem of the form in (16). Among numerical methods, the choice is small due to complicated properties of the functions $F_{ij}$.

At the same time, the functions $F_{ij}$, $i \in \{\mathrm{I}, \mathrm{II.a}, \mathrm{II.b}\}$, $j \in \{1,2\}$ satisfy the so-called the $\varepsilon$-Lipschitz condition (see [4]) in the corresponding set $W_i$. Now we give the corresponding definition.

Let $f$ be a function from a subset $A$ contained in a normed space $X$ to a normed space $Y$. The function $f$ is called *$\varepsilon$-Lipschitz continuous* on $A$, if for any number $0 < \varepsilon < \varepsilon_0$ there exists some number $L(\varepsilon) \geqslant 0$ such that for all $x, y \in A$, the following inequality holds

$$\| f(x) - f(y) \| \leqslant L(\varepsilon) \| x - y \| + \varepsilon. \tag{17}$$

By $\varepsilon$-$\mathrm{Lip}(A; Y)$ denote the set of $\varepsilon$-Lipschitz continuous on $A$ functions; for a fixed function $f$

and fixed $\varepsilon$, by $lip(f;\varepsilon)$ denote the smallest function (i.e., the infimum) of all $L(\varepsilon)$ satisfying (17). (The properties for $lip(f;\varepsilon)$ can be found in [2] and [4], where it is denoted by $l(\varepsilon)$) and $\kappa(\varepsilon)$, respectively).

To present the assertions that each function $F_{ij}$, $i\in\{\mathrm{I},\mathrm{II.a},\mathrm{II.b}\}$, $j\in\{1,2\}$ is $\varepsilon$-Lipschitz continuous on a certain set in $\mathbb{R}^n$ ($n$ equal to 2 or 3), we shall give the explicit formulas for $\mathrm{Pr}_1 h_{kj}^{-1}(N)$, $\mathrm{Pr}_2 h_{kj}^{-1}(N)$, $W_k$, $k\in\{\mathrm{I},\mathrm{II}\}$, $j\in\{1,2\}$.

Using (3) and (7), we get

$$\varphi(N):=\mathrm{Pr}_1 h_{kj}^{-1}(N)=\arctan(y/x), k\in\{\mathrm{I},\mathrm{II}\}, j\in\{1,2\};$$

$$\psi_{\mathrm{I}\,j}(N):=\mathrm{Pr}_2 h_{\mathrm{I}j}^{-1}(N)=\frac{\pi}{4}+(-1)^j\left(\frac{\pi}{4}-\frac{1}{2}\arcsin\frac{\sqrt{x^2+y^2}}{v^2}\right), j\in\{1,2\};$$

$$\begin{aligned}\psi_{\mathrm{II}\,j}(N):\ &=\mathrm{Pr}_2 h_{\mathrm{II}j}^{-1}(N)=\\ &=\arctan\left((x^2+y^2)^{-1/2}\left(v^2+(-1)^j\sqrt{v^4-(x^2+y^2+2v^2z)}\right)\right), j\in\{1,2\}.\end{aligned}\tag{18}$$

Note there exist some values of the variables $x$, $y$, and $z$ such that the expression under the last square root in (18) can be negative. But this fact indicates that, for the given magnitude $v_0$ of the initial projectile velocity, the coordinates of the projectile drop point can not be equal to those $x$, $y$, and $z$ for any values of $\psi$.

The $\varepsilon$-Lipschitz continuity of the functions $F_{ij}$, $i\in\{\mathrm{I},\mathrm{II.a},\mathrm{II.b}\}$, $j\in\{1,2\}$, is formulated below by Propositions 1–3. We also give several lemmas in which some functions that are used in the definitions of the functions $F_{ij}$ are $\varepsilon$-Lipschitz continuous or Lipschitz continuous. For the subsequent estimates of $lip(\cdot)$ and $lip(\cdot;\varepsilon)$, we use the norm $\|\cdot\|_1$; we also assume that $g_1, g_2\in \mathrm{Lip}([\theta_1,\theta_2];\mathbb{R})$, $H\in\mathrm{Lip}(\mathbb{R}^3;\mathbb{R})$, $2\varepsilon\in(0,\pi/2]$.

**Lemma 1.** *For the function $g(\varphi,\psi)$ defined by (6), we have $g\in\mathrm{Lip}(\Theta;\mathbb{R})$ and*

$$lip(g)\leqslant\max\{lip(g_1), lip(g_2), 1\}.\tag{19}$$

The proofs of this lemma and the following ones can be found above in the Russian section.

**Lemma 2.** *For each $k\in\{\mathrm{I},\mathrm{II}\}$, it is true that $\varphi(x,y)\in\mathrm{Lip}(W_k;(-\pi/2,\pi/2))$ and*

$$lip(\varphi)\leqslant\sqrt{2}\kappa^{-1}.\tag{20}$$

**Lemma 3.** *For each $j\in\{1,2\}$, it is true that $\psi_{I\,j}(x,y)\in\varepsilon\text{-}\mathrm{Lip}(W_{\mathrm{I}};(0,\pi/2))$ and*

$$lip(\psi_{\mathrm{I}\,j};\varepsilon)\leqslant\begin{cases}\left(\sqrt{2}v^2\sqrt{1-\tau^2(2\varepsilon)}\right)^{-1} & \text{if}\quad 0<2\varepsilon<\pi/2-1,\\ \left(\sqrt{2}v^2\right)^{-1}(\pi/2-2\varepsilon) & \text{if}\quad \pi/2-1\leqslant 2\varepsilon\leqslant\pi/2,\end{cases}\tag{21}$$

*where $\tau(2\varepsilon)$ is the unique root of the equation $(\pi/2-2\varepsilon-\arcsin\tau)\sqrt{1-\tau^2}=1-\tau$ in the interval $[0,1)$.*

**Proposition 1.** *For each $j\in\{1,2\}$, it is true that $F_{Ij}(N)\in\varepsilon\text{-}\mathrm{Lip}(W_I;\mathbb{R})$ and*

$$lip(F_{\mathrm{I}j};\varepsilon)\leqslant\max\left\{\omega_{\mathrm{I}j1}lip(g)\left(lip(\varphi)+lip\left(\psi_{\mathrm{I}\,j};\varepsilon(\omega_{\mathrm{I}j1}lip(g))^{-1}\right)\right),\omega_{\mathrm{I}j2}\right\},\tag{22}$$

*where $lip(g)$, $lip(\varphi)$, and $lip\left(\psi_{\mathrm{I}\,j};\varepsilon(\omega_{\mathrm{I}j1}lip(g))^{-1}\right)$ determined by (19), (20), and (21), respectively.*

**Lemma 4.** *For each $j\in\{1,2\}$, it is true that $\psi_{\mathrm{II}j}(N)\in\varepsilon\text{-}\mathrm{Lip}(W_{\mathrm{II}};(-\pi/2,\pi/2))$ and*

$$lip\left(\psi_{\mathrm{II}j};\varepsilon\right)\leqslant\left(\rho_{\mathrm{II}}(2\varepsilon)^{-1}+\beta\sqrt{2}\right)\kappa^{-2},\tag{23}$$

*where* $\beta = v^2 + \sqrt{v^4 - \kappa^2 - 2v^2 z_{\min}}$.

**Proposition 2.** *For each* $j \in \{1,2\}$, *it is true that* $F_{\mathrm{II.a}\,j}(N) \in \varepsilon\text{-}\mathrm{Lip}(W_{\mathrm{II}}; \mathbb{R})$ *and*

$$lip(F_{\mathrm{II.a}\,j}; \varepsilon) \leqslant \max\left\{\omega_{\mathrm{II.a}\,j1} lip(g)\left(lip(\varphi) + lip\big(\psi_{\mathrm{II}j}; \varepsilon(\omega_{\mathrm{II.a}\,j1} lip(g))^{-1}\big)\right), \omega_{\mathrm{II.a}\,j2}, \omega_{\mathrm{II.a}\,j3}\right\}, \quad (24)$$

*where* $lip(g)$, $lip(\varphi)$, *and* $lip\big(\psi_{\mathrm{II}j}; \varepsilon(\omega_{\mathrm{II.a}\,j1} lip(g))^{-1}\big)$ *determined by (19), (20), and (23), respectively.*

**Proposition 3.** *For each* $j \in \{1,2\}$, *it is true that* $F_{\mathrm{II}.b\,j}(N) \in \varepsilon\text{-}\mathrm{Lip}(W_{\mathrm{II}}; \mathbb{R})$ *and*

$$\begin{aligned} lip(F_{\mathrm{II.b}\,j}; \varepsilon) \leqslant \max \quad &\Big\{\omega_{\mathrm{II.b}\,j1} lip(g)\left(lip(\varphi) + lip\big(\psi_{\mathrm{II}\,j}; \varepsilon(\omega_{\mathrm{II.b}\,j1} lip(g))^{-1}\big)\right), \omega_{\mathrm{II.b}\,j2}, \\ &\omega_{\mathrm{II.b}\,j3} lip(H), \omega_{\mathrm{II.b}\,j4} lip(H), \omega_{\mathrm{II.b}\,j5} lip(H)(\rho_{\mathrm{II}} lip(H)/(8\varepsilon) + 1)\Big\}, \end{aligned} \quad (25)$$

*where* $lip(g)$, $lip(\varphi)$, *and* $lip\big(\psi_{\mathrm{II}\,j}; \varepsilon(\omega_{\mathrm{II.b}\,j1} lip(g))^{-1}\big)$ *determined by (19), (20), and (23), respectively.*

According to our assumptions, we have $M \in W_i$, $i \in \{\mathrm{I}, \mathrm{II}\}$, whence $F_{ij}(M) \geqslant 0$, $j \in \{1,2\}$. These inequalities and Propositions 1–3 imply that problem (16) of finding a projection can be solved by using the algorithms proposed in [2, 3]. Each algorithm either generates an infinite or a finite sequence $Q^m$. In the first case, the sequence $Q^m$ converges to a zero of function $F_{ij}$. This zero is closest to $M$ under the norm, and the sequence also satisfies inequality $F_{ij}(Q^m) > 0$ for each $m = 1,2,\ldots$. In the second case, the last element of the sequence is the required point $N_*$, in which $F_{ij}(N_*) = 0$. In this paper we use one of these algorithms. Before describing it, as we said above, we turn to the weighting coefficients $\omega_{ijl}$ used in constructing the functions $F_{ij}$, because choosing $\omega_{ijl}$ affects the accuracy of the solution obtained by the algorithm.

Since the algorithm finds the next approximation to the optimal solution by using information about possible change in the value of the function $F_{ij}$ with respect to the current approximation and since the function $F_{ij}$ is defined by maximum of some function-components, we must take into account how these function-components are related to each other. If these relations are ignored, then the approximations generated by the algorithm may not reach the stopping criterion $F_{ij}(Q^m) < \varepsilon^*$ within a reasonable time (see also remarks 1–2 below). There is a reason for it: the accuracy for a function-component of $F_{ij}$ may be too high, whence the number of algorithm iterations to reach accuracy is too large. It is known (see [5]) that an upper bound of the number of iterations for the gradient-free algorithms, including our algorithm, is determined by using possible change in the value of the objective function; therefore, if we know these changes for all function-components of $F_{ij}$, then we can take $\omega_{ijl}$ to provide the number of iterations adequate to the overall accuracy $\varepsilon^*$. In other words, having taken a value of the overall accuracy $\varepsilon^*$ for the values of $F_{ij}$, we can specify accuracies for every function-component of $F_{ij}$ by using the values of the corresponding weighting coefficients: accuracy for $l$-th function-component will be equal to $\varepsilon^*/\omega_{ijl}$.

In this article, we suppose $\omega_{ij1} = 1$, $i \in \{\mathrm{I}, \mathrm{II.a}, \mathrm{II.b}\}$, $j \in \{1,2\}$, that is, the overall accuracy $\varepsilon^*$ corresponds to the accuracy of the first function-component of $F_{ij}$. The accuracies $\omega_{ijl}$ for other function-components are taken according to remark above.

Let us present the algorithm.

**Algorithm** for solving Problems I, II.a, and II.b

*Step 0.* Choose: an initial value $\varepsilon_0 > 0$ of the $\varepsilon$-Lipschitz parameter, the initial projectile velocity magnitude determined by $v$ or $v_0$, a lower bound $\kappa$ of $x$, an initial point $N^0 = M \in W_{ij}$, parameters $\gamma, \lambda \in (0,1)$. Assign the visibility zone constraints determined by $g(\varphi, \psi)$ satisfying (6). Set $Q^0 := N^0$, $k := 0$, and $m := 0$.

*Step 1.* Depending on Problem I, II.a, or II.b, calculate $F_{ij}(N^k)$, $i \in \{\mathrm{I}, \mathrm{II.a}, \mathrm{II.b}\}$, $j \in \{1,2\}$ using (12), (13), or (15), respectively. If $F_{ij}(N^k) < \varepsilon_k(1+\gamma)$, then sequentially set $N^{k+1} := N^k$, $Q^{m+1} := N^k$, and $m := m+1$ and pass to Step 2. Otherwise, go to Step 3.

*Step 2.* If $F_{ij}(N^k) \leqslant \varepsilon_k$, then set $\varepsilon_{k+1} := \lambda F_{ij}(N^k)$. Otherwise, set $\varepsilon_{k+1} := \lambda\varepsilon_k$. In either case, go to Step 4.

*Step 3.* Find $N^{k+1}$ by the following scheme. Depending on Problem I, II.a, or II.b, calculate $lip(F_{ij}, \varepsilon_k)$, $i \in \{\mathrm{I}, \mathrm{II.a}, \mathrm{II.b}\}$, $j \in \{1,2\}$ by formulas (22), (24), or (25), respectively, using (19)–(21), and (23). Solve the problem

$$N^{k+1} = \operatorname{argmin}\{F_{ij}(X) | X \in \mathrm{fr}K_k \cap W_{ij}\},\ \ i \in \{\mathrm{I}, \mathrm{II.a}, \mathrm{II.b}\}, j \in \{1,2\}, \qquad (26)$$

where

$$K_k = \{X \in \mathbb{R}^n |\ \| X - M \| \leqslant \textstyle\sum_{s=0}^{k} r_s\}, r_k = \frac{F_{ij}(N^k) - \varepsilon_k}{\sqrt{n}\, lip(F_{ij}; \varepsilon_k)},$$

$n = 2$ for Problem $i = \mathrm{I}$ and $n = 3$ for Problems $i = \mathrm{II.a}$ and $i = \mathrm{II.b}$. Next, set $\varepsilon_{k+1} := \varepsilon_k$ and go to Step 4.

*Step 4.* If $F_{ij}(N^{k+1}) = 0$ or $F_{ij}(Q^m) = 0$, then $N^{k+1}$ or $Q^m$ is regarded as a solution of problem (16), respectively. Otherwise, set $k := k+1$ and pass to Step 1.

**Remark 1.** According to the main statement about the algorithm convergence (see [2]), if the number of points $N^k$ is infinity, then the algorithm ensures the convergence $F_{ij}(Q^m) \xrightarrow[m\to\infty]{} 0$. Therefore, the condition of the form $F_{ij}(Q^m) < \varepsilon^*$, $\varepsilon^* \in (0,1)$, should be included in the stopping criterion. We can also include the conditions setting accuracies for some function-components of $F_{ij}$ or the condition of the form $\| Q^m - Q^{m+1} \| < \varepsilon_Q^*$, $\varepsilon_Q^* \in (0,1)$ in the stopping criterion.

**Remark 2.** It is clear that the main computational costs are produced by Step 3, where the optimization subproblem for finding the point $N^{k+1}$ is solved. We will give two pieces of advice about solving this subproblem. First, since each function $F_{ij}$, $j \in \{1,2\}$, is minimized either on a circle (Problem $i = \mathrm{I}$) or on a 2-sphere (Problems $i = \mathrm{II.a}, \mathrm{II.b}$), it follows that the dimension of the subproblem can be reduced by one by applying transformation to polar coordinates or spherical ones, respectively. Secondly, it is not necessary to solve this subproblem precisely. Indeed, since the value $F_{ij}(N^k)$ is compared with $\varepsilon_k$ and $\varepsilon_k(1+\gamma)$, instead of the point $N^{k+1}$ determined by (26), one can find a point $\widetilde{N}^{k+1}$ such that $F_{ij}(N^{k+1}) \leqslant F_{ij}(\widetilde{N}^{k+1}) < F_{ij}(N^{k+1}) + \delta_k$, where $\delta_k < \varepsilon_k\gamma$.

**Remark 3.** For computing values of the 4-th function-component of $F_{\mathrm{II}.b}$ in problem II.b, we must solve the one-dimensional optimization problem of minimization of the function $f(\lambda) = H(\lambda N + (1-\lambda)M)$. This problem can be solved some corresponding numerical method. A similar remark also concerns the 5-th function-component of $F_{\mathrm{II}.b}$.

## 3. Numerical Examples

The algorithm was tested on several examples in which the target $M$, the set $W_i$, $i \in \{\mathrm{I}, \mathrm{II}\}$, and the overall accuracy $\varepsilon^*$ that taken equal to an initial value $\varepsilon_0$ of the $\varepsilon$-Lipschitz parameter were varied. The optimization subproblem for finding the point $N^{k+1}$ (see Step 3) was solved by the uniform search method.

For all examples, the following parameter values were set: $v_0 = 180\text{m/s}$; $\kappa = 100\text{m}$; $\gamma = \lambda = 0.5$. The values of the variable parameters are presented in Tables 1–3. Hereinafter and in Tables 1–3, we use the following notation: $M_1 = (110,0,20)$, $M_2 = (2700,0,-10)$ are the target $M$ (the third coordinate of $M$ for Problem I was equal to 0); $E_1 = \{(\varphi,\psi) | \varphi \in [0, \pi/2 - 10^{-3}], \psi \in [7\pi/36, 8\pi/36]\}$, $E_2 = \{(\varphi,\psi) | \varphi \in [0, \pi/2 - 10^{-3}], (4 + \sin\varphi)\pi/36 \leqslant \psi \leqslant (1 + \sin\varphi)\pi/9\}$; $(x_N, y_N, z_N)$ is the solution of Problem I, II.a, or II.b; $(\varphi,\psi)$ is the angles for the direction of the gun barrel, which ones corresponding to $(x_N, y_N, z_N)$; $k_{\text{tot}}$ is the total number of algorithm iterations; $t$ is the total computational time in seconds. The value $z_{\min}$ in Problems II.a and II.b was equal to – 10, and the set $D$ in Problems II.b was determined by function

$$H(x,y,z) = \min\{\max\{90 - x, x - 130, -10 - y, y - 30, z - 20\}, z + 10\}.$$

The weighting coefficients in the functions $F_{\text{I}}$, $F_{\text{II.a}\,j}$, and $F_{\text{II.b}\,j}$ (see (12), (13), and (15)) were chosen as follows: I) $\omega_{\text{I}j1} = 1$, $\omega_{\text{I}j2} = 0.01$; II.a) $\omega_{\text{II.a}\,j1} = 1$, $\omega_{\text{II.a}\,j2} = \omega_{\text{II.a}\,j3} = 0.01$; II.b) $\omega_{\text{II.b}\,j1} = 1$, $\omega_{\text{II.b}\,j2} = 0.01$, $\omega_{\text{II.b}\,j3} = \omega_{\text{II.b}\,j4} = \omega_{\text{II.b}\,j5} = 0.001$.

The tests were performed on an IntelCorei3-4020 1.50 GHz personal computer. The numerical results are presented in Tables 1–3. It follows that the algorithm finds the $\varepsilon$-solution in Problem I faster than in other two Problems. Moreover, in Problem II.a and especially in Problem II.b, there is a large increase in the number of iterations with $\varepsilon_0$ greater than 0.05. (The dash in Table 3 means that the algorithm did not give a solution for 200sec.) A reasonable explanation for this increase may be that, for the functions $F_{\text{II.a\,,b}}$ a decrease in $\varepsilon$ implies an increase in $L(\varepsilon)$. However, for a finite number, number of iterations (when $L(\varepsilon)$ is not very high), it is still possible to cut off a set that does not contain zero of $F_{ij}$; this yields that the set $W_i$ is narrowed down to its subset. It is readily seen that computational time can be reduced by parallelizing some parts of the algorithm. For example, one can be applied in subproblems, in which the uniform search method are used.

| $M$ | $E$ | $\varepsilon$ | $x_N$, m | $y_N$, m | $\varphi$, ° | $\psi$, ° | $k_{\text{tot}}$ | $t$, sec |
|---|---|---|---|---|---|---|---|---|
| $M_1$ | $E_1$ | 0.1 | 2818.6 | 0.0 | 0.0 | 29.3 | 300 | 0.146 |
| | | 0.05 | 2976.8 | 0.0 | 0.0 | 32.1 | 353 | 0.170 |
| | | 0.01 | 3081.4 | 0.0 | 0.0 | 34.4 | 575 | 0.306 |
| | $E_2$ | 0.1 | 1523.4 | -152.1 | -5.7 | 13.8 | 301 | 0.419 |
| | | 0.05 | 1834.6 | -91.4 | -2.9 | 16.9 | 442 | 0.822 |
| | | 0.01 | 2068.1 | -20.6 | -0.6 | 19.4 | 891 | 3.514 |
| $M_2$ | $E_1$ | 0.1 | 2820.3 | 0.0 | 0.0 | 29.3 | 40 | 0.034 |
| | | 0.05 | 2977.1 | 0.0 | 0.0 | 32.2 | 144 | 0.094 |
| | | 0.01 | 3081.4 | 0.0 | 0.0 | 34.4 | 331 | 0.213 |
| | $E_2$ | 0.1 | 2606.4 | 48.4 | 1.1 | 26.0 | 36 | 0.031 |
| | | 0.05 | 2447.8 | 138.6 | 3.2 | 24.0 | 158 | 0.086 |
| | | 0.01 | 2316.7 | 217.0 | 5.4 | 22.4 | 456 | 0.244 |

Table 1: Results of test examples for Problem I.

| $M$ | $E$ | $\varepsilon$ | $x_N$, m | $y_N$, m | $z_N$, m | $\varphi$, ° | $\psi$, ° | $k_{\text{tot}}$ | $t$, sec |
|---|---|---|---|---|---|---|---|---|---|
| $M_1$ | $E_1$ | 0.1 | 100.0 | 0.1 | 58.0 | 0.1 | 31.0 | 1889 | 1.249 |
| | | 0.05 | 100.0 | 1.8 | 65.4 | 1.0 | 34.1 | 4182 | 2.244 |
| | | 0.01 | 102.0 | 2.1 | 72.1 | 1.2 | 36.1 | 26297 | 5.227 |
| | $E_2$ | 0.1 | 108.7 | -0.5 | 25.9 | -0.3 | 14.3 | 614 | 1.073 |
| | | 0.05 | 106.8 | -1.0 | 31.4 | -0.6 | 17.3 | 2444 | 1.495 |
| | | 0.01 | 104.8 | -0.3 | 35.9 | -0.2 | 19.8 | 17619 | 8.854 |
| $M_2$ | $E_1$ | 0.1 | 2730.7 | 0.0 | 48.3 | 0.0 | 29.3 | 8478 | 4.064 |
| | | 0.05 | 2771.2 | 0.0 | 121.4 | 0.0 | 32.2 | 33069 | 12.750 |
| | | 0.01 | 2803.2 | -0.1 | 173.5 | 0.0 | 34.4 | 213343 | 103.261 |
| | $E_2$ | 0.1 | 2502.2 | 107.1 | -8.6 | 2.5 | 24.4 | 27281 | 5.952 |
| | | 0.05 | 2336.3 | 826.1 | -4.3 | 19,5 | 24.2 | 244288 | 131.227 |
| | | 0.01 | 2324.5 | 663.2 | -2.8 | 15,9 | 23.4 | 313626 | 163.866 |

Table 2: Results of test examples for Problem II.a

| $M$ | $E$ | $\varepsilon$ | $x_N$, m | $y_N$, m | $z_N$, m | $\varphi$, ° | $\psi$, ° | $k_{\text{tot}}$ | $t$, sec |
|---|---|---|---|---|---|---|---|---|---|
| $M_1$ | $E_1$ | 0.1 | 100.0 | 0.0 | 58.0 | 0.0 | 31.0 | 2748 | 9.753 |
| | | 0.05 | 100.2 | 16 | 62.6 | 9.1 | 32.6 | 7308 | 14.240 |
| | $E_2$ | 0.1 | 108.7 | -0.5 | 25.9 | -0.3 | 14.3 | 1072 | 11.303 |
| | | 0.05 | 106.8 | -1.0 | 31.4 | -0.6 | 17.3 | 4258 | 11.179 |
| $M_2$ | $E_1$ | 0.1 | 2730.7 | 0.0 | 48.3 | 0.0 | 29.3 | 12746 | 19.626 |
| | | 0,05 | - | - | - | - | - | - | - |
| | $E_2$ | 0.1 | 2502.2 | 107.1 | -8.6 | 2.5 | 24.4 | 47995 | 88.068 |
| | | 0.05 | 2437.0 | 145.7 | -10.0 | 3.4 | 23.6 | 69230 | 121.723 |

Table 3: Results of test examples for Problem II.b

Kazan National Research Technical University named after A.N. Tupolev-KAI (Russia).
*E-mail address*: NKArutyunova@kai.ru, AMDulliev@kai.ru, VIZabotin@kai.ru